\documentclass{amsart}[11pt]
\usepackage{amsmath,amscd}
\usepackage[utf8]{inputenc}
\usepackage{multicol}
\usepackage{url}
\usepackage{xcolor}
\usepackage{colortbl}
\usepackage{multirow}
\usepackage{subcaption}

\usepackage{footnote}
\makesavenoteenv{tabular}
\makesavenoteenv{table}

\usepackage{algorithm}
\usepackage{algpseudocode}

\usepackage{booktabs, array}
\usepackage{enumitem}
\usepackage{amsthm}
\usepackage{mathtools}
\usepackage{comment}
\theoremstyle{plain}
\newtheorem{theorem}{Theorem}[section]
\newtheorem{lemma}[theorem]{Lemma}
\newtheorem{proposition}[theorem]{Proposition}

\theoremstyle{definition}

\newtheorem{example}[theorem]{Example}

\newtheorem{definition}[theorem]{Definition}
\newtheorem{remark}[theorem]{Remark}

\newtheorem{hypothesis}[theorem]{Hypothesis}

\numberwithin{theorem}{section}
\numberwithin{equation}{section}
\numberwithin{table}{section}
\numberwithin{figure}{section}
\numberwithin{algorithm}{section}

\usepackage{url}
\usepackage[colorlinks=true,citecolor=blue,linkcolor=blue,urlcolor=blue,bookmarks,bookmarksopen,bookmarksdepth=2,backref=page,breaklinks]{hyperref}


\renewcommand*{\backref}[1]{}
\renewcommand*{\backrefalt}[4]{%
	\ifcase #1 Not cited.%
	\or        p.~#2.%
	\else      pp.~#2.%
	\fi}

\def\F{\mathbb F}

\def\Z{\mathbb Z}

\def\rmc(#1,#2){RM(#1,#2)}
\def\form(#1,#2){H(#1,#2)}
\def\val(#1,#2){V(#1,#2)}
\def\boole(#1){ B(#1) }
\def\bst(#1,#2,#3){\boole(#1,#2,#3)}
\def\tildeboole(#1){ \widetilde B(#1) }
\DeclareMathOperator{\comp}{Comp}
\def\qfspace(#1){{{\rm Q}(#1)}}
\def\flat{{\mathfrak f}}
\def\tfrbis(#1,#2){ {#1}^\dag(#2) }
\def\fd{{\mathbb F}_2}

\def\Flat(#1){{{\mathfrak F}(#1)}}
\def\supp(#1){{\rm supp}(#1)}

\def\der(#1,#2){{\rm Der}_{#2} (#1)}

\def\stab(#1){\textsc{stab}(#1)}
\def\stablevel(#1,#2){\textsc{stab}_{#1}(#2)}

\def\rmq(#1,#2,#3){\rmc(#1,#2)/\rmc(#3,#2)}

\def\walsh(#1,#2){\widehat{#1}(#2)}

\def\level(#1){\underset{#1}{=}}

\def\binogauss(#1,#2){\genfrac{[}{]}{0pt}{}{#1}{#2}}

\def\val(#1){{\rm val}(#1)}

\def\anf(#1){{\rm anf}(#1)}
\def\fix(#1,#2,#3,#4){{{\rm fix}^{#1,#2}_{#3}}{(#4)}}
\def\classe(#1,#2,#3){{{\mathcal T}(#1,#2,#3)}}

\def\stab(#1){\textsc{stab}(#1)}
\def\stablevel(#1,#2){\textsc{stab}^{#1}(#2)}
\def\stableveldim(#1,#2,#3){\textsc{stab}^{#1}_{#2}(#3)}
\def\stableveldeg(#1,#2,#3,#4){\textsc{stab}^{#1,#2}_{#3}(#4)}
\def\level(#1){\underset{#1}{\sim}}
\def\bound(#1,#2){\underset{#1}{\overset{#2}{\sim}}}
\def\modulo(#1,#2){\mod\rmc(#1,#2)}

\def\pow#1.#2{\tiny$10^{#1.#2}$}

\def\level(#1){\underset{#1}{=}}
\def\val(#1){{\rm val}(#1)}

\def\anf(#1){{\rm anf}(#1)}
\def\fix(#1,#2,#3,#4){{{\rm fix}^{#1,#2}_{#3}}{(#4)}}
\def\classe(#1,#2,#3){{{\mathcal T}(#1,#2,#3)}}

\def\stab(#1){\textsc{stab}(#1)}
\def\stablevel(#1,#2){\textsc{stab}^{#1}(#2)}
\def\stableveldim(#1,#2,#3){\textsc{stab}^{#1}_{#2}(#3)}
\def\stableveldeg(#1,#2,#3,#4){\textsc{stab}^{#1,#2}_{#3}(#4)}
\def\level(#1){\underset{#1}{\sim}}
\def\bound(#1,#2){\underset{#1}{\overset{#2}{\sim}}}
\def\modulo(#1,#2){\mod\rmc(#1,#2)}

\def\pow#1.#2{\tiny$10^{#1.#2}$}
\def\degV(#1,#2,#3){\deg_{#1+#2}(#3)}

\DeclareMathOperator*{\argmax}{argmax}
\def\ortho(#1){{#1}^{\star}}

\title[Millions of inequivalent APN functions]{Millions of inequivalent quadratic APN functions in eight variables}

\author[C. Beierle]{Christof Beierle}
\author[P. Langevin]{Philippe Langevin}
\author[G. Leander]{Gregor Leander}
\author[A. Polujan]{Alexandr Polujan}
\author[S. Rasoolzadeh]{Shahram Rasoolzadeh}

\address{Ruhr University Bochum, Germany}
\email{christof.beierle@rub.de}
\email{gregor.leander@rub.de}
\email{shahram.rasoolzadeh@rub.de}

\address{Imath, University of Toulon, France}
\email{philippe.langevin@univ-tln.fr}

\address{Otto-von-Guericke-Universit\"{a}t Magdeburg, Germany}
\email{alexandr.polujan@gmail.com}

\date{\today}

\begin{document}

\begin{abstract}
The only known example of an almost perfect nonlinear (APN) permutation in even dimension was obtained by applying CCZ-equivalence to a specific quadratic APN function in dimension six. Motivated by this result, there have been numerous attempts to construct new quadratic APN functions. Prior to this work, $32\,892$ quadratic APN functions in dimension eight are known and two recent conjectures address their possible total number. The first, proposed by Y. Yu and L. Perrin, suggests that there are more than $50\,000$ such functions. The second, by A. Polujan and A. Pott, argues that their number exceeds that of inequivalent quadratic $(8,4)$-bent functions, which is $92\,515$. We computationally construct $3\,775\,599$  inequivalent quadratic APN functions in dimension eight and estimate the total number to be about six million.
\end{abstract}

\maketitle
\section{Introduction}

Let $\F_2^n$ be the $n$-dimensional vector space over the finite field with two elements $\F_2$. A mapping $F\colon \F_2^n \to \F_2^m$ is referred to as an \emph{$(n,m)$-function}. The \emph{differential uniformity} of an $(n,m)$-function $F$ is defined as
\begin{equation}\label{def:diff_uniformity}
	\delta_F = \max_{\substack{a \in \mathbb{F}_2^n \setminus \{0\} \\ b \in \mathbb{F}_2^m}} \left|\{x \in \mathbb{F}_2^n \colon F(x+a) + F(x) = b\}\right|.
\end{equation}
An $(n,m)$-function $F$ is said to be \emph{almost perfect nonlinear (APN)} if it achieves the lowest possible differential uniformity, that is, if $\delta_F = 2$. Note that for any $a \neq 0$ and any $b$, the cardinality of the set $\{x \in \F_2^n \colon F(x+a) + F(x) = b\}$ is always even, since solutions occur in pairs $x$ and $x+a$.

Originally introduced in the $(n,n)$-case in the context of resistance against differential attacks~\cite{DBLP:conf/crypto/NybergK92}, APN functions play an important role in cryptography. Of particular interest are APN permutations on $\F_2^n$ in even dimensions $n$, whose existence was posed as an open problem in~\cite{Dobbertin1999:Niho}. Since then, it has been shown that APN permutations do not exist in dimension four~\cite{Hou06}, while the only known example has been constructed in dimension six~\cite{Dillon-perm}. The existence of such permutations in higher even dimensions is known as the \textsc{BIG APN problem}.  Moreover, APN functions (including those that are permutations) have important connections to coding theory~\cite{DBLP:journals/dcc/CarletCZ98}, incidence structures in finite geometry~\cite{Yoshiara10}, difference sets~\cite{Dillon1999,DillonDobbertin2004}, and Sidon sets~\cite{CarletPicek2021,CzerwinskiP26} in combinatorics. APN functions in the $(n,m)$-case with $n<m$, again defined by the property $\delta_F=2$, have recently attracted increased attention due to their use in secondary constructions of $(n,n)$-APN functions~\cite{DBLP:journals/dcc/BeierleLP22,DBLP:journals/corr/abs-2501-03922}. For additional results on APN functions, we refer to~\cite{DBLP:journals/corr/abs-2501-03922,Carlet2021_Book}.

The investigation of construction methods of APN functions began with the study of infinite families of \emph{monomial functions} $x \mapsto x^d$ over the finite fields $\mathbb{F}_{2^n}$~\cite{Gold1968,JanwaWilson1993}. Only six infinite families of APN monomials are currently known (see Table~\ref{table:apn_monomials}), and a complete classification remains open even within this restricted class. While several notions of equivalence are used for this purpose, \emph{Carlet-Charpin-Zinoviev (CCZ)} equivalence~\cite{DBLP:journals/dcc/CarletCZ98} is the most general one commonly considered in the study of APN functions. Two $(n,m)$-functions $F$ and $F'$ are said to be CCZ-equivalent if there exists an affine permutation $\mathcal{L}$ on $\F_2^{n} \times \F_2^{m}$ such that $\mathcal{L}\left(\mathcal{G}_{F}\right)=\mathcal{G}_{F'}$, where $\mathcal{G}_{F}=\{(x,F(x))\colon x\in\F_2^n\}$ denotes the \emph{graph} of $F$.

\begin{table}
	\centering
	\caption{Known infinite families of APN power functions $x^d$ on $\mathbb{F}_{2^n}$}
	\begin{tabular}{c|cc|c}
		\hline
		Family & Exponent $d$ & Conditions & References \\
		\hline 
		Gold & $2^i + 1$ & $\gcd(i,n)=1$ & \cite{Gold1968,Nyberg93} \\
		\hline
		Kasami & $2^{2i} - 2^i + 1$ & $\gcd(i,n)=1$ & \cite{JanwaWilson1993, Kasami1971} \\
		\hline
		Welch & $2^t + 3$ & $n = 2t + 1$ & \cite{Dobbertin1999:Welch} \\
		\hline
		\multirow{2}{*}{Niho} 
		& $2^t + 2^{t/2} - 1,\ t\ \text{even}$ & \multirow{2}{*}{$n = 2t + 1$} & \multirow{2}{*}{\cite{Dobbertin1999:Niho}} \\
		\cline{2-2}
		& $2^t + 2^{(3t+1)/2} - 1,\ t\ \text{odd}$ &  &  \\
		\hline
		Inverse & $2^{n} - 2$ & $n = 2t + 1$ & \cite{BethDing1994,Nyberg93} \\
		\hline
		Dobbertin & $2^{4i} + 2^{3i} + 2^{2i} + 2^i - 1$ & $n = 5i$ & \cite{Dobbertin2001:DivBy5} \\
		\hline
	\end{tabular}
	\label{table:apn_monomials}
\end{table}
	
The question of whether all APN functions are CCZ-equivalent to monomials was raised in~\cite{DobbertinMMPW2003}. Shortly thereafter, the first examples of quadratic APN functions that are CCZ-inequivalent to monomials were discovered~\cite{EdelKP2006}. These developments further stimulated the investigation of new quadratic APN functions, especially since two important examples of APN functions were obtained by modifying quadratic ones. Among them are the APN permutation in dimension six~\cite{Dillon-perm}, obtained by applying a CCZ-equivalence to a quadratic function, and an APN function that is neither CCZ-equivalent to a quadratic function nor to a monomial~\cite{EdelP09}, constructed via \emph{switching} from a quadratic function, that is, by modifying a single component function. Since then, research on quadratic APN functions has been largely driven by the following two challenges:
\begin{itemize}
	\item[(I)] discovering new examples and constructions, and 
	\item[(II)] resolving equivalence problems among the known functions.
\end{itemize}

Over the past 15 years, significant progress has been made in addressing these questions. Several infinite families of APN functions have been constructed; see, e.g.,~\cite{GologluK2025_CT,Li_Kaleyski_2024} for a detailed overview. Some of these families have been shown to be CCZ-inequivalent to permutations~\cite{GologluL20,GologluPavlu2021}. Notably, one family of quadratic APN functions~\cite{DBLP:journals/dcc/Taniguchi19a} was shown to contain exponentially many CCZ-inequivalent functions as $n$ grows~\cite{KaspersZhou2021}. At the same time, numerous \emph{sporadic examples}, that is, functions that do not (yet) belong to any known infinite family, have been discovered through computational search techniques. The majority of these functions were identified in~\cite{DBLP:journals/corr/abs-2009-07204,DBLP:journals/dcc/BeierleLP22,DBLP:journals/dcc/YuWL14}. These discoveries have, in turn, driven the development of more efficient methods for constructing and classifying APN functions~\cite{DBLP:journals/ccds/Kaleyski21,CanteautCP2022}. In Table~\ref{table:apn_table_plot}, we present some key results on the computational construction and classification of APN functions. 

\begin{table}
	\centering
	\caption{Known classification results for APN functions on $\mathbb{F}_2^n$}
	\begin{tabular}{c|l|c}
		\hline
		$n$                      & APN functions (\# of CCZ-classes)                                                                 & References                                                                                                                                \\ \hline
		$n \leq 5$               & All APN functions                                                           & \cite{DBLP:journals/dcc/BrinkmannL08}                                                                                    \\ \hline
		$n = 6$                  & All quadratic ($13$) and all cubic ($1$)  & \cite{DBLP:journals/dcc/Pott16} and \cite{langevin}                                                             \\ \hline
		$n = 7$                  & All quadratic ($488$)                             & \cite{DBLP:journals/dcc/YuWL14,DBLP:journals/iacr/KalginI20a}                                                                                     \\ \hline
		\multirow{2}{*}{$8\le n\le 11$} & \multirow{2}{*}{Quadratic\footnote{Note that the results for $n \geq 8$ are not exhaustive.} ($\approx 33\,000$)} & \cite{DBLP:journals/ffa/BeierleCLP22,DBLP:journals/corr/abs-2009-07204,DBLP:journals/dcc/BeierleLP22,EdelP09,IchanskaBKY2025} \\ 
		&                                                                             & \cite{DBLP:journals/dcc/Taniguchi19a,weng2013quadratic,YuPerrin2024,DBLP:journals/dcc/YuWL14}\\ \hline
		$n=8$ & Quadratic ($3\,775\,599$) & \textit{This article} \\ \hline
	\end{tabular}
	\label{table:apn_table_plot}
\end{table}

Up to dimension $5$, the search space is small enough to allow for a complete classification. From dimension $6$ onward, additional constraints, such as low algebraic degree, are required. These constraints have enabled the classification of all quadratic and cubic APN functions in dimension $6$, as well as all quadratic ones in dimension $7$, and have led to the discovery of numerous sporadic examples in higher dimensions. A conjecture in~\cite{YuPerrin2024}, based on the fraction of the explored search space, suggested that the number of CCZ-inequivalent APN functions in dimension $8$ exceeds $50\,000$. In this paper, we present two methods for generating a large number of quadratic APN functions, which is far exceeding what was previously known or predicted. In total, using substantial computational resources, we generated more than $3.5$ million inequivalent APN functions.

The two approaches to generating quadratic APN functions that we use follow a common principle: starting from partial functions and extending them step by step to larger ones. However, the ways in which this is done are fully ``orthogonal''.  In the first, and by far more successful, method, we generate functions by extending an $(8,m)$-function to an $(8,m+1)$-function, i.e., by adding a coordinate function. The key starting point is a vectorial $(8,4)$-function that is \emph{bent}, whose output space is then extended one dimension at a time. This approach is motivated by the recent classification of all $(8,4)$-bent functions given in~\cite{PolujanPott2022BFA}. The second approach proceeds by starting with an $(n,8)$-function that is APN and extending it to an $(n+1,8)$-APN function, i.e., by increasing the input space one dimension at a time. The advantage of this approach is that it is arguably more random and does not rely on the hypothesis that any quadratic APN function arises as an extension of an $(8,4)$-bent mapping. The downside is that it is even more computationally demanding; consequently, only a small fraction of the millions of newly constructed APN functions were generated in this way.

The paper is organized as follows. In Section~\ref{sec:preliminaries}, we provide the necessary background on vectorial Boolean functions, including APN and bent functions. In Section~\ref{sec:millions}, we consider an approach for constructing millions of APN functions by extending vectorial bent functions. In particular, in Section~\ref{sec:hypothesis}, we introduce the working hypothesis that every quadratic APN function on $\mathbb{F}_2^8$ can be obtained as an extension of an $(8,4)$-bent function, and we provide evidence for its correctness in this dimension based on results in small dimensions, known classifications, and probabilistic arguments. In Section~\ref{sec: bent to apn}, we describe the \emph{``coordinate approach''} used to construct millions of APN functions by extending quadratic vectorial bent functions, and we discuss the resulting APN functions. In Section~\ref{sec: randomized apn}, we introduce the recursive \emph{``variable approach''} based on increasing the number of variables in APN functions while preserving the APN property; it is worth noting that with this approach, we classified all quadratic $(n,8)$-APN functions, for $1\le n\le 6$. In Section~\ref{sec:results}, we consider the constructed APN functions in more detail. In particular, Section~\ref{sec:constructed_apns} discusses several important properties of the newly constructed APN functions in comparison with the previously known ones: CCZ-equivalence to an APN permutation, possible distributions of component functions, and the orders of the automorphism groups. Then, in Section~\ref{sec:total_number}, using the functions constructed via the ``variable approach,'' we provide a theoretical estimate of the total number of inequivalent APN functions in dimension eight, which is approximately six million.

\section{Preliminaries}\label{sec:preliminaries}
Let $n,m$ be positive integers. An \emph{$(n,m)$-function} $F$ is a mapping from $\mathbb{F}_2^n$ into $\mathbb{F}_2^m$, in particular, $(n,1)$-functions are referred to as \emph{Boolean functions}, while $(n,m)$-functions with $m\ge2$ are referred to as \emph{vectorial (Boolean) functions}. For an $(n,m)$-function $F$, we define for each $b \in \mathbb{F}_2^m$ a Boolean function $F_b \colon x \mapsto  b \cdot F(x)$, which is called a \emph{component function} of $F$; here ``$\cdot$'' denotes the usual dot product on $\mathbb{F}_2^m$. We define the \emph{space of components} of $F$ as the $\mathbb{F}_2$-vector space $$\operatorname{Comp}(F) = \{ F_b \colon b \in \mathbb{F}_2^m \}.$$
For a Boolean function $f$ on $\F_2^n$, the \emph{Walsh transform} of $f$ at the point $a\in\F_2^n$ is defined by 
$$\quad\hat{\chi}_{f}(a)=\sum_{x\in\F_2^n} (-1)^{f(x)+ a\cdot x} \in \Z.$$
The \emph{Walsh transform} of a vectorial $(n,m)$-function $F$ is defined as the mapping \[\hat{\chi}_{F}\colon \F_2^n\times\F_2^m\rightarrow \Z,\quad (a,b) \mapsto \hat{\chi}_{F_b}(a).\] Every Boolean  function $f$ on $\F_2^n$ has a unique representation as a multivariate polynomial
\begin{equation*}\label{ANF}
	f(x_1, x_2, \ldots, x_n ) = \sum_{S\subseteq \{1,2,\ldots, n\}} a_S X_S,
	\quad  \text{with } \ a_S\in\fd, \ X_S = \prod_{s\in S} x_s,
\end{equation*}
which is called the \emph{algebraic normal form (ANF)}. The \emph{(algebraic) degree} of a non-zero function $f$, denoted by $\deg(f)$, is the maximal cardinality of $S$ with  $a_S=1$ in the ANF of $f$. We have $\deg(0) = - \infty$ by convention. We say that $f$ is \emph{homogeneous} if all the terms in its ANF have the same degree. Every vectorial $(n,m)$-function $F$ can be written as $F(x) = (f_1(x), \dots, f_m(x))$, for all $x\in\F_2^n$, where each Boolean function $f_i$ on $\F_2^n$ is called a \emph{coordinate function}. In turn, the ANF of a vectorial $(n,m)$-function $F$ is defined coordinate-wise. Consequently, the degree of $F$ is the maximum degree among its coordinate functions, i.e., $\deg(F)=\max_{1\le i\le m} \deg(f_i)$. Functions of degree at most one are called \emph{affine}, and those of degree two are \emph{quadratic}.

In this article, we consider two important classes of $(n,m)$-functions possessing minimal differential uniformity: (vectorial) bent functions and APN functions, which are defined as follows.

\begin{definition}\label{def:apn}
	Let $n \le m$. An $(n,m)$-function $F$ is called \emph{almost perfect nonlinear (APN)} if its differential uniformity satisfies $\delta(F)=2$.
\end{definition}

\begin{definition}\label{def:bent}
	An $(n,m)$-function $F$ is called \emph{bent} (or \emph{perfect nonlinear}) if its differential uniformity satisfies $\delta(F)=2^{n-m}$.
\end{definition}

\begin{remark}
	Boolean bent functions $f\colon \F_2^n \to \F_2$ exist only for even $n$, see~\cite{ROTHAUS1976}. These functions are characterized by the property that the Walsh transform of $f$ satisfies
	$\hat{\chi}_f(a)=\pm 2^{n/2}$, for all $a\in\F_2^n$. Vectorial $(n,m)$-bent functions $F$ are characterized by the property that, for every non-zero $b\in\F_2^m$, the component function $F_b$ is Boolean bent. Equivalently, the Walsh transform of $F$ satisfies
	$\hat{\chi}_F(a,b)=\pm 2^{n/2}$ for all $a\in\F_2^n$ and all $b\in\F_2^m\setminus\{0\}$. Vectorial $(n,m)$-bent functions exist only if $m \le n/2$; this condition is referred to as the \emph{Nyberg bound}~\cite{Nyberg91}.
\end{remark}

In order to get a characterization of APN functions using the Walsh transform, we use the following quantity. For a Boolean function $f$ on $\F_2^n$, the \emph{normalized fourth power moment of the Walsh transform}, denoted $\alpha(f)$, is defined by
\begin{equation}\label{eq: alpha}
    \alpha(f) = \frac{1}{2^{3n}}\sum_{u \in \mathbb{F}_2^n} (\hat{\chi}_f(u))^4.
\end{equation}

\begin{remark}
   1.  Note that a Boolean function $f$ is bent on $\F_2^n$ if and only if $\alpha(f) = 1$. Consequently, a vectorial $(n,m)$-function $F$ is bent if and only if 
    \begin{equation}\label{eq:bent_alpha}
    	\sum_{0\neq f \in \operatorname{Comp}(F) } \alpha(f) = 2^m - 1.
    \end{equation} 
    
   \noindent 2. A quadratic function $f$ on $\F_2^n$ given by $f(x)=xUx^T+l(x)$, where $U$ is an upper triangular $n\times n$-matrix with zero diagonal and $l$ is an affine function on $\F_2^n$,  satisfies $$\alpha(f)=2^{n-\operatorname{rank}(f)},$$ where $\operatorname{rank}(f)=\operatorname{rank}_{\F_2}(U+U^T)$. In the case $n=8$, which we are interested in, the correspondence between the normalized fourth power moment of the Walsh transform and the rank of a quadratic Boolean function (which is always even~\cite[Chapter 15]{MacWilliamsSloane}) is as follows:
    $$\alpha(s) = 
		\begin{cases}
			1 &\mbox{if } \operatorname{rank}(s)=8 \\
			4 &\mbox{if } \operatorname{rank}(s)=6 \\
			16 &\mbox{if } \operatorname{rank}(s)=4 \\
			64 &\mbox{if } \operatorname{rank}(s)=2
		\end{cases}.$$
\end{remark}

Similarly to Boolean and vectorial bent functions, APN functions on $\F_2^n$ can be characterized using the normalized fourth power moment of the Walsh transform of their component functions in the following way.

\begin{theorem}\cite[Corollary 1]{berger2006almost}
	An $(n, n)$-function $F$ is APN if and only if 
	\begin{equation}\label{eq:APN_alpha}
		\sum_{0\neq f \in \operatorname{Comp}(F) } \alpha(f) = 2^{n+1} - 2.
	\end{equation}
\end{theorem}	

Besides the notion of \emph{CCZ-equivalence} given in the previous section, another important notion is extended-affine equivalence. We say that $(n,m)$-functions $F$ and $F'$ are \emph{extended-affine equivalent} (\emph{EA-equivalent}) if there exist affine permutations $A$ and $B$ of $\F_2^m$ and $\F_2^n$, respectively, and an affine $(n,m)$-function $C$ such that $F' = A \circ F \circ B + C$. While, in general, CCZ-equivalence is a coarser equivalence relation than EA-equivalence, for quadratic $(n,n)$-APN functions, all (vectorial) bent functions, and Boolean functions, these two notions coincide; for results on the corresponding classes of functions, we refer to~\cite{yoshiara2008dimensional},~\cite{KyureghyanP08}, and~\cite{BudaghyanC2010}, respectively.

Motivated by the partitioning of quadratic APN functions into EA-equivalence classes, the following notion was introduced in~\cite{CanteautCP2022}.
For a quadratic APN $(n,n)$-function $F$ on $\F_2^n$, there exists a unique function $\pi_F$ on $\F_2^n$ such that $\pi_F(0) = 0$, $\pi_F(a) \neq 0$ for $a \neq 0$, and for any $(a, x) \in (\mathbb{F}_2^n)^2$, it holds that
$$
\pi_F(a) \cdot (F(x + a) + F(x) + F(a) + F(0) ) = 0.
$$
Such a function $\pi_F$ is called the \emph{ortho-derivative} of $F$. Crucially, the EA-equivalence of two quadratic APN $(n,n)$-functions implies the EA-equivalence of their ortho-derivatives~\cite[Proposition 36]{CanteautCP2022}. Using invariants for the equivalence on the ortho-derivative, one can efficiently deduce the  EA-inequivalence of (even large numbers of) quadratic APN functions, as indicated in~\cite{DBLP:journals/corr/abs-2009-07204,YuPerrin2024}.

\section{Millions of quadratic APN functions via the coordinate approach}\label{sec:millions}
The coordinate approach is a well-known construction method of vectorial functions with prescribed properties by successively adding new coordinate functions to a given function. For vectorial functions whose defining properties can be characterized through the structure of their component space, this approach has been highly successful in many cases, for example in the construction and classification of vectorial bent functions~\cite{Arshad_PhD_Thesis,PolujanPott2022BFA,PolujanP21}, whose component space forms a vector space of Boolean bent functions. For quadratic APN functions, however, this approach is less straightforward, partly because the structure of the component space is not well understood in general.

In this section, we adapt the coordinate approach to construct millions of inequivalent quadratic APN functions in dimension eight. Our approach is based on the assumption that the component space of any quadratic APN function in dimension eight contains a bent subspace of dimension four, i.e., that such functions arise as extensions of vectorial $(8,4)$-bent functions. We then develop an efficient algorithm for constructing these extensions.

\subsection{Working hypothesis on the existence of a bent subspace in the component space}\label{sec:hypothesis}

We focus on constructing $(n, n)$-APN functions by extending known vectorial $(n, m)$-functions with $m < n$, through the addition of new coordinate functions. To formalize this construction process, we introduce the notion of extension in the following way.

\begin{definition}
	An $(n,r)$-function $G$ is called an \emph{extension} of an $(n,m)$-function $F$ if the component space $\operatorname{Comp}(F)$ is a subspace of $\operatorname{Comp}(G)$.
\end{definition}

For $n=4$, there exists a unique quadratic $(4,2)$-bent function up to EA-equivalence~\cite{DingMT19}, which can be extended to the unique quadratic APN function $x \mapsto x^3$ over $\mathbb{F}_{2^4}$~\cite{DBLP:journals/dcc/BrinkmannL08}. For $n=6$, there are three EA-inequivalent quadratic $(6,3)$-bent functions~\cite[Theorem 3.1]{Arshad_PhD_Thesis}, each admitting multiple extensions to functions among the 13 EA-inequivalent quadratic APN functions on $\mathbb{F}_2^6$. In particular, each of these 13 functions is an extension of a $(6,3)$-bent function, as shown in~\cite[Table 3.1]{Arshad_PhD_Thesis}. These observations motivate the following working hypothesis in dimension eight:

\begin{hypothesis}\label{hypothesis}
    A quadratic $(8,8)$-APN function is an extension of an $(8,4)$-bent function. 
\end{hypothesis}

To challenge this hypothesis, we confirmed that the $32\,892$ known quadratic APN functions in dimension 8 (prior to this work; see the list and references in~\cite[Sec.\@ 4]{DBLP:journals/dcc/BeierleLP22}) contain a bent space of dimension 4. To do so, there are algorithms to check for the existence of large vector spaces contained in a given set~\cite{DBLP:conf/asiacrypt/BonnetainPT19,GologluPavlu2021}.

Besides that confirmation, the strongest argument for the correctness of the hypothesis stems from the fact that any quadratic APN function in an even number of variables $n$ always consist of at least $\frac{2}{3}(2^n-1)$ bent components (i.e., $170$ bent components for $n=8)$, see~\cite[Section 4.2]{DBLP:conf/fse/Nyberg94}. If we take a set $A$  sampled uniformly at random from the set of all subsets of $\mathbb{F}_2^8 \setminus \{0\}$ of size 170, the existence of a $4$-dimensional vector space contained in $A \cup \{0\}$ is almost guaranteed. More precisely, if we denote by $N$ the random variable corresponding to the number of $k$-dimensional subspaces contained in $A \cup \{0\}$, where $A$ is a uniformly chosen $t$-subset of $\mathbb{F}_2^n \setminus \{0\}$, we obtain
\[\mathbb{E}(N) = \binom{n}{k}_2 \cdot \frac{\binom{t}{2^k-1}}{\binom{2^n-1}{2^k-1}}.\]
Here, $\binom{n}{k}_2$ denotes the Gaussian binomial coefficient, which equals the number of $k$-dimensional subspaces contained in $\mathbb{F}_2^n$. For the parameter $$(n,k,t) = (8,4,170),$$ this yields an expected value of $\mathbb{E}(N) \approx 369.38$. Adapting the method of Clark and Suen~\cite{clark1999probability}, it is also possible to precisely evaluate $\mathbb{E}(N^2)$ and to derive a theoretical lower bound on $\mathbb{P}(N \geq 1)$ by the second moment method. For our case, we obtain  \[\mathbb{P}(N \geq 1) \geq 0.978.\]

However, this lower bound seems far from tight as our experiments suggest: We sampled subsets $A$ of $\F_2^8 \setminus \{0\}$ of size 170 uniformly at random and checked the existence of a 4-dimensional vector space  within $A \cup \{0\}$. In $10^6$ experiments, we always detected such a subspace.

\begin{remark}
    We have no reason to expect that Hypothesis~\ref{hypothesis} holds for even $n$ in general (in the sense that a quadratic APN function in  $n$ variables is an extension of an $(n,\frac{n}{2})$-bent function). Indeed, for the following cases $(n,k,t)=(n,\frac{n}{2},\frac{2}{3}(2^n-1))$ with
    \begin{equation*}
    	\begin{split}
    		(n,k,t)& = (10,5,682), \\
    		(n,k,t)&= (12,6,2\,730), \\
    		(n,k,t)&= (14,7,10\,922), \\
    	\end{split}
    \end{equation*}
     we have $\mathbb{E}(N) \approx 300.78$, $\mathbb{E}(N) \approx 1.46$, and $\mathbb{E}(N) < 10^{-7}$, respectively. For $n = 10$, the second moment method still yields $\mathbb{P}(N \geq 1) \geq 0.971$, while for $n=12$ and $n=14$, we can only deduce $\mathbb{P}(N \geq 1) \geq 0.58$ and $\mathbb{P}(N \geq 1) \geq 6 \cdot 10^{-8}$, respectively.
\end{remark}

\subsection{Extending quadratic vectorial bent functions}\label{sec: bent to apn}
Let $\qfspace(n)$ denote the space of quadratic functions on $\mathbb{F}_2^n$. Quadratic vectorial $(n,m)$-bent functions are precisely the $m$-dimensional subspaces of $\qfspace(n)$ such that every non-zero element has rank $n$. Recently, all $92\,515$ quadratic $(8,4)$-bent functions were classified in~\cite{PolujanPott2022BFA} by successively extending $(8,2)$-bent functions to $(8,3)$-bent functions, and then to $(8,4)$-bent functions. However, further extension to $(8,8)$-APN functions using the ``bent-template'' is impossible due to the Nyberg bound~\cite{Nyberg91}. That is, for each additional coordinate function, only a subset of all component functions will be bent. To capture these differences, we introduce the notion of the profile of a vectorial function.

\begin{definition}
    For a subspace $V\subset\qfspace(n)$, define the pair $\operatorname{BS}(V) = (b, K(V))$ where $b$ is the number of bent functions in $V$ and $$K(V) = \sum_{0 \neq s \in V} \alpha(s).$$ The \emph{profile} $\operatorname{P}_k(F)$ of a quadratic $(n,m)$-function $F$ is the lexicographically greatest tuple
    \begin{equation}\label{eq: profile}
        \operatorname{P}_k(F)=[\operatorname{BS}(V_0), \operatorname{BS}(V_1), \operatorname{BS}(V_2), \ldots, \operatorname{BS}(V_k)],
    \end{equation}
	where $(V_0 \subset V_1 \subset \cdots \subset V_k)$ ranges over flags of $\operatorname{Comp}(F)$ satisfying $\dim(V_i)=i$.
\end{definition}

\begin{remark}
    Given a monotonic function $w$ for the inclusion of subspaces in a vector space of dimension $m$, the goal is to determine a \textit{maximal flag} in the lexicographical order induced by $w$. To achieve this, we use a recursive approach combined with  \textit{branch pruning}. The algorithm explores possible subspaces by adding linearly independent vectors while maintaining a table $F$, where each $F[k]$ represents the weight of the subspace of dimension $k$. It also maintains a table $\operatorname{opt}$, where $\operatorname{opt}[k]$ denotes the $k$-th entry of the best profile found so far. At each step, the search is pruned if $F[k] < \operatorname{opt}[k]$, and the optimal flag is updated if $F[k] > \operatorname{opt}[k]$. The recursive procedure continues until all possible dimensions are explored, ensuring the discovery of the maximal flag in the lexicographical order.
\end{remark}

\begin{remark} It turns out that the $32\,892$ known quadratic APN functions (prior to this work) share relatively few distinct $P_6$ profiles, which we describe in Table~\ref{table: profiles}. As pointed out in the previous section, each of those functions can be derived from a quadratic $(8,4)$-bent function, i.e., their $P_4$ profiles are equal to $[(0, 0), (1, 1), (3, 3), (7, 7), (15, 15)]$.
\end{remark}
\begin{table}[h!]
	\centering
	\caption{Profiles $P_6(F)$ of $32\,892$ known quadratic $(8,8)$-APN functions $F$}
	\begin{tabular}{@{}rc@{}} \toprule
		$\#$APN & Profile Sequence \\ \midrule
		9       & [(0, 0), (1, 1), (3, 3), (7, 7), (15, 15), (28, 40), (50,           102)] \\
		468     & [(0, 0), (1, 1), (3, 3), (7, 7), (15, 15), (28, 40), (52, \phantom{1}96)] \\
		90      & [(0, 0), (1, 1), (3, 3), (7, 7), (15, 15), (28, 40), (54, \phantom{1}90)] \\
		32      & [(0, 0), (1, 1), (3, 3), (7, 7), (15, 15), (30, 34), (50,           102)] \\
		3,128   & [(0, 0), (1, 1), (3, 3), (7, 7), (15, 15), (30, 34), (52, \phantom{1}96)] \\ 
		7       & [(0, 0), (1, 1), (3, 3), (7, 7), (15, 15), (30, 34), (54,           102)] \\
		13,480  & [(0, 0), (1, 1), (3, 3), (7, 7), (15, 15), (30, 34), (54, \phantom{1}90)] \\ 
		2       & [(0, 0), (1, 1), (3, 3), (7, 7), (15, 15), (30, 34), (56,           108)] \\
		7,549   & [(0, 0), (1, 1), (3, 3), (7, 7), (15, 15), (30, 34), (56, \phantom{1}84)] \\
		49      & [(0, 0), (1, 1), (3, 3), (7, 7), (15, 15), (30, 34), (56, \phantom{1}96)] \\ 
		7,008   & [(0, 0), (1, 1), (3, 3), (7, 7), (15, 15), (30, 34), (58, \phantom{1}78)] \\
		80      & [(0, 0), (1, 1), (3, 3), (7, 7), (15, 15), (30, 34), (58, \phantom{1}90)] \\ 
		923     & [(0, 0), (1, 1), (3, 3), (7, 7), (15, 15), (30, 34), (60, \phantom{1}72)] \\
		67      & [(0, 0), (1, 1), (3, 3), (7, 7), (15, 15), (30, 34), (60, \phantom{1}84)] \\ \bottomrule
	\end{tabular}
    \label{table: profiles}
\end{table}

\begin{example}
    Let $F\colon\F_{2^8}\to\F_{2^8}$ be defined by  $x \mapsto x^3$. The profile $P_8$ of $F$ is given by
    $$P_8(F)=[(0, 0), (1, 1), (3, 3), (7, 7), (15, 15), (28, 40), (50, 102), (90, 238), (170, 510)].$$
\end{example}

\begin{remark}
    We observe that $x^3$ is indeed an extension of an $(8,4)$-bent function, and that its profile requires $28$ bent functions in dimension 5. More generally, the profile of $x^3$ is closely related to the distribution of non-cubes in subspaces. In dimension $2t$, the Boolean function $$x\in\F_{2^{2t}} \mapsto \operatorname{Tr}(b x^3)$$ is a bent component of $x^3$ if and only if $b$ is not a cube in $\F_{2^{2t}}$. In particular, when $t$ is odd, $x^3$ is an extension of a $(2t,t)$-bent function; see~\cite[Section 3]{Carlet11}. When $t$ is even, there are no general results, to the best of our knowledge, except for some estimates about the distribution of cubes/non-cubes considered in the paper~\cite{GologluL20}. This raises an important question about the normality of the indicator of cubes: does there exist an $\F_2$-subspace of dimension $t$ of $\F_{2^{2t}}$ that contains no non-zero cubes?
\end{remark}

At certain steps, we need to distinguish between different extensions using invariants for EA-equivalence, since full classification is too time- and resource-consuming due to a large number of functions involved. The following invariant has proven to be highly efficient in terms of both efficient discrimination and computational cost.

\begin{definition}
     Let ${\rm Q}_2(n)\subset \qfspace(n)$ be the set of all quadratic homogeneous functions of rank 2 in $n$ variables. Let $F$ be a quadratic $(n,m)$-function. Define the multiset
	\begin{equation}
		J_2(F) = \{* \; \operatorname{rank}(q + f) \colon q \in {\rm Q}_2(n), \; f \in \operatorname{Comp}(F) \; *\},
	\end{equation}
	where $\operatorname{Comp}(F) = \{ F_b \colon b \in \mathbb{F}_2^m \}$ is the space of component functions of $F$.
\end{definition}

We show that $J_2$ is an invariant under EA-equivalence for quadratic $(n,m)$-functions.

\begin{proposition}
	Let $F$ and $F'$ be two EA-equivalent quadratic $(n,m)$-functions. Then, $J_2(F) = J_2(F')$.
\end{proposition}
\begin{proof}
	Since $F$ and $F'$ are EA-equivalent, there exist an invertible $m \times m$ matrix $A$, an invertible $n \times n$ matrix $B$, a vector $c \in \mathbb{F}_2^n$, and an affine function $C \colon \mathbb{F}_2^n \to \mathbb{F}_2^m$ such that
		$$F'(x) = A(F(xB + c)) + C(x).$$
	For any $b \in \mathbb{F}_2^m$, the component function $f'(x) = b \cdot F'(x)$ is given by:
		$$f'(x) = b \cdot [A(F(xB + c)) + C(x)] =  f(xB+c) + C_b(x),$$
		where $f=F_{A^Tb}\in\operatorname{Comp}(F)$.
	The term $C_b(x)$ is affine and does not affect the rank of the quadratic part. Similarly, translating the input by $c$ contributes only an affine term, leaving the quadratic part unchanged. Hence, $f'(x)$ and $f(xB)$ have identical quadratic parts.
	
	Now, consider the multiset $J_2(F')$. Its elements are the ranks of functions of the form $q'(x) + f'(x)$ where $q' \in {\rm Q}_2(n)$ and $f' \in \operatorname{Comp}(F')$. Using the change of variables $x = yB^{-1}$, we obtain
	$$\operatorname{rank}(q'(x) + f'(x)) = \operatorname{rank}(q'(yB^{-1}) + f(y))=\operatorname{rank}(q(y) + f(y)),$$
	where $q(y) = q'(yB^{-1})$. Since the map $q' \mapsto q$ induces a permutation of ${\rm Q}_2(n)$ and $f' \mapsto f$ is bijective from $\operatorname{Comp}(F')$ to $\operatorname{Comp}(F)$, it follows that  $J_2(F) = J_2(F')$.
\end{proof}

\begin{remark}
    For the set of quadratic $(8,5)$-functions having $28$ bent components, we identified $3\,747\,371\,328$ distinct $J_2$ labels. This gives a lower bound on the number of EA-equivalence classes of such functions.
\end{remark}

\begin{lemma}
    \label{DIFF}
    An $(n,m)$-function $F$ having differential uniformity greater than $2^{n-m+1}$ does not have any APN extension. 
\end{lemma}
\begin{proof}
It is well-known that an $(n,m-1)$-function obtained by deleting a coordinate of an $(n,m)$-function with differential uniformity $\delta$ has differential uniformity at most $2\delta$. Applying this process $n-m$ times to an $(n,n)$-APN function (for which the differential uniformity is $\delta=2$), any $(n,m)$-function $F$ obtained in this way satisfies $\delta_F\le 2^{\,n-m+1}$.
\end{proof}

In order to introduce an efficient extension procedure for constructing APN functions, it is convenient to express the APN property through linear-algebraic conditions on the component space of a function, which leads to the following definition.

\begin{definition}
	Let $V$ be a subspace of the space of Boolean functions in $n$ variables equipped with a scalar product $(f,g) \mapsto \langle f,g \rangle$. We say that $v\in V$ is a \emph{flat} vector if  there exists four distinct  $x$, $y$, $z$, $t$ in $\mathbb{F}_2^n$ such that $x+y+z+t = 0$ and 
	$$\forall  f\in V, \quad \langle v, f \rangle = f(x)+f(y)+f(z)+f(t).$$
	We denote by $\Flat(V)$ the set of flat vectors of $V$ and by $\flat$ the indicator of this set.
\end{definition}	

We now introduce a characterization of the APN property based on flat vectors in spaces of Boolean functions.
    
\begin{lemma} \label{FLAT}
	A vectorial function $F \colon \mathbb{F}_2^n \rightarrow \mathbb{F}_2^m$ whose components are in $V$ is APN if and only if $ \comp(F)^\perp \cap \Flat(V)=\emptyset$.
\end{lemma}

\begin{proof}
    Let $S$ be the component space of $F$.
    By definition, $F$ is APN if and only if, for every 4-subset $\{x, y, z, t\}$ of $\mathbb{F}_2^n$ with $x+y+z+t = 0$, the sum $\sigma = F(x)+F(y)+F(z)+F(t)$ is non-zero in $\mathbb{F}_2^m$. Equivalently, there exist a component function $f \in S$ such that $$f(x)+f(y)+f(z)+f(t) = 1.$$ By the definition of a flat vector $v$ representing the set $\{x, y, z, t\}$, this requirement is exactly $\langle v, f \rangle = 1$. Thus, $F$ is APN if and only if every $v \in \Flat(V)$ satisfies $v \notin S^\perp$. This condition is equivalent to $$ S^\perp \cap \Flat(V)  = \emptyset,$$ which completes the proof.
\end{proof}

	Let $S= \operatorname{Comp}(F)$ and $S^\perp$ be its orthogonal complement. 
	For the indicator $\flat$ of the set $\Flat(V)$,  its Fourier transform is given by $$\widehat{\flat}(s) = \sum_{v \in V} \flat(v) (-1)^{\langle v, s \rangle }.$$ 
	Applying the Poisson's summation formula to the indicator $\flat$ of the set $\Flat(V)$, we get
	\begin{equation*}
		\begin{split}
			\frac{1}{|S|} \sum_{s \in S} \widehat{\flat}(s) &=\frac{1}{|S|} \sum_{s \in S} \sum_{v \in V} \flat(v) (-1)^{{\langle v, s \rangle }}\\
            &= \sum_{v \in V} \flat(v) \left( \frac{1}{|S|} \sum_{s \in S} (-1)^{{\langle v, s \rangle }} \right)\\
			&=\sum_{v \in V} \flat(v) 1_{S^\perp}(v) = \sum_{r \in S^\perp} \flat(r).
		\end{split}
	\end{equation*}
	Since $\flat$ is the indicator of $\Flat(V)$, the sum $\sum_{r \in S^\perp} \flat(r)$ simply counts how many elements of $S^\perp$ are also in $\Flat(V)$. In particular, 
    $F$ is APN if and only if the following condition holds
    \begin{equation}\label{eq: zero}
	   \sum_{s\in  \operatorname{Comp}(F) 
       } \widehat{\flat}(s) = 0.
    \end{equation}

\begin{remark}
	In the case where $V$ is the space of homogeneous quadratic functions in $n$ variables, we 
    consider the scalar product with respect to the basis $\{x_i x_j \colon 1 \leq i < j \leq n\}$. By definition,
	$$ x_ix_j .  \left( \sum_{r, s} a_{rs}\, x_rx_s\right) = a_{ij}.$$ By the M\"obius inversion formula, the coefficient $a_S$ in the ANF of a Boolean function $f = \sum a_S X_S$ is given by $$a_S = \sum_{\supp(v) \subseteq S} f(v).$$ Specifically, for a homogeneous quadratic function $f$, the scalar product with the basis element $x_i x_j$ yields
	\[ \langle x_i x_j, f \rangle = a_{ij}=f(0) + f(e_i) + f(e_j) + f(e_i + e_j), \]
	since for $S=\{i,j\}$, its subsets correspond to the four vectors $0,e_i,e_j,e_i+e_j$. This identifies $x_i x_j$ as a flat vector associated with the linear subspace $$\{0, e_i, e_j, e_i + e_j\}.$$ 
    
    More generally, every quadratic form of rank $2$ is EA-equivalent to some $x_ix_j$, meaning it can be written as $x_ix_j \circ A$ for some invertible linear map $A$. Under such a change of variables, the associated 4-subset transforms to $$\{0, e_iA^{-1}, e_jA^{-1}, (e_i+e_j)A^{-1}\},$$ which still sums to zero, so the flatness property is preserved.
    One can see that
    two affine planes having the same direction define the same flat 
    character, hence, by a counting argument,
    the set $\Flat(V)$ corresponds to the set of quadratic 
    forms of rank 2, which are the forms EA-equivalent to $x_i x_j$. 
    
    Under this identification, Lemma~\ref{FLAT} recovers the known result of Edel~\cite[Lemma 8]{Edel2011} that a quadratic function $F$ is APN if and only if $\comp(F)^\perp$ contains no quadratic form of rank 2.
\end{remark}

Lemma~\ref{FLAT} and  Relation \ref{eq: zero}, can be used to determine all the APN extensions $E$ of a given $(n, n-r)$ function $F$  with $\operatorname{Comp}(E) \subseteq V$. To do this, we must choose a subspace $T$ of dimension $r$ in $W:= V/\operatorname{Comp}(F)$ such that 
\begin{equation}\label{EXTENSION}
	\sum_{t \in T}  \sigma(t)= 0,
    \quad
    \text{where}
    \ 
    \sigma(t) := \sum_{s\in\operatorname{Comp}(F)} \widehat{\flat}(t+s).
\end{equation}

\begin{remark}\label{rem:is_extendable}
	In the case $r=2$, finding $T$ satisfying the condition in Eq.~\eqref{EXTENSION} is feasible when the dimension of $W$ is relatively small (e.g., up to 30). We have to find $u,v\in W$ that form a basis of $T$ such that:
	\[
	\sigma(0) + \sigma(u) + \sigma(v) + \sigma(u + v) = 0.
	\]
	We may assume without loss of generality that $\sigma(u) \leq \sigma(v) \leq \sigma(u+v)$, which implies
	\begin{equation}\label{eq:extension_condition1}
	    \sigma(u) \leq -\frac{1}{3} \sigma(0).
	\end{equation}
    After finding such $u$, we need to find $v$ satisfying the condition
    \begin{equation}\label{eq:extension_condition2}
        2\sigma(v) \leq \sigma(v) + \sigma(u+v) = -\sigma(0) - \sigma(u).
    \end{equation}
\end{remark}

While the conditions above provide an efficient test for determining whether a given function is APN-extendable, the following observation is useful for constructing APN extensions directly via linear algebra.

\begin{remark}\label{rem:construct}
    Let $W$ be any \emph{complement} of $\operatorname{Comp}(F)$ in $V$, i.e., $V=\operatorname{Comp}(F)\oplus W$. Note that $\operatorname{Comp}(F)^\perp$ is not necessarily a complement of $\operatorname{Comp}(F)$ in $V$. Assume that an $(n,n-1)$-function $F$ is extendable to an $(n,n)$-APN function. To construct an $(n,n)$-APN extension, one can solve the linear system over $\mathbb{F}_2$ arising from Lemma~\ref{FLAT} and Relation~\ref{eq: zero}. To do so, we enumerate all 2-dimensional flats $\{ x,y,z,t \}$ of $F$ s.t. $$F(x)+F(y)+F(z)+F(t)=0$$ and impose the constraint $$w(x) + w(y) + w(z) + w(t) = 1$$ on the unknown coordinate $w\in W$. This technique is referred to as \emph{switching}~\cite{EdelP09}. Assuming that $V$ is a space of quadratic homogeneous functions in $n$ variables and $\operatorname{Comp}(F)\subset V$, we have that $$\dim(W)=\binom{n}{2} - (n-1) = \binom{n-1}{2}.$$ The resulting linear system over $\F_2$ in these $\binom{n-1}{2}$ free variables of $W$ is solved by Gaussian elimination; a unique solution (corresponding to a quadratic APN function) exists if and only if the system has full rank $\binom{n-1}{2}$ and is consistent.
\end{remark}

\begin{remark}
Given a quadratic APN mapping, one can use a similar approach to investigate the existence of non-quadratic switchings. All switchings of the quadratic APN mappings constructed in this work are themselves quadratic, in contrast to the case $n=6$, where switchings of degree $3$ do exist.
\end{remark}

Using the ideas above, we can describe the extension procedure of quadratic vectorial bent functions used to obtain quadratic APN functions having a given profile 
\begin{equation}
\label{TARGET}
[(0, 0), (1, 1), (3, 3), (7, 7), (15, 15), (30,34), (B_6, S_6)],\qquad \text{with}\quad B_6 \geq 54 
\end{equation}
as Algorithm~\ref{alg:APNfromBent}. Using this method, quadratic APN extensions of a quadratic $(8,6)$-function can be computed in approximately $0.05$ seconds on a usual desktop computer.

\begin{algorithm}
	\caption{Quadratic $(8,8)$-APN functions from $(8,4)$-bent functions}
	\begin{algorithmic}[1]
		\Require The set $\mathcal{B}$ of all quadratic $(8,4)$-bent functions~\cite{PolujanPott2022BFA}.
		\Ensure A set of EA-inequivalent $(8,8)$-APN functions satisfying Profile \ref{TARGET}.
		
		\State\textbf{Step 1:} For each bent function in $\mathcal{B}$ compute all $(8,5)$-extensions with $B_5=30$ bent functions.
		
		\State\textbf{Step 2:} Apply the $J_2$-invariant to the obtained extensions. This yields a set containing $2\,403\,534$ functions.
		
		\State \textbf{Step 3:} Compute all $(8,6)$-extensions with $B_6 \ge 54$ bent functions, which have the Profile \ref{TARGET}. 
		
		\For{each obtained extension}
		\State (1) Apply the $J_2$ invariant to the resulting set of $(8,6)$-functions to reduce
        \State the search space.
		\State (2) Use Lemma~\ref{DIFF} to eliminate functions with too large differential
        \State uniformity.
		\State (3) For each $v \in W = \comp(F)^\perp$, accept $v$ as a seventh coordinate function
        \State (yielding an $(8,7)$-extension) only if the condition in Eq.~\eqref{eq:extension_condition1} can still be 
        \State satisfied for some $u \in W$ satisfying the condition in Eq.~\eqref{eq:extension_condition2} (Remark~\ref{rem:is_extendable}).
        \State (4) Apply the $J_2$ invariant to the resulting set of $(8,7)$-functions to reduce
        \State the search space.
        \State (5) Use Remark~\ref{rem:construct} to construct $(8,8)$-APN functions.
		\EndFor
		
		\State\textbf{Step 4:} Apply the ortho-derivative to the obtained $(8,8)$-APN functions to establish pairwise CCZ-inequivalence.
		
		\State \Return Constructed $(8,8)$-APN functions.
	\end{algorithmic}
	\label{alg:APNfromBent}
\end{algorithm}

\begin{remark}
    The procedure described in Algorithm~\ref{alg:APNfromBent} is efficient when $B_6 > 54$, in the case of $B_6 = 54$, we extend only half of the $(8,6)$-functions. Yet, the direction of the profile specified by Eq.~\eqref{TARGET} provides almost all the APN functions we found.
\end{remark}

\section{Thousands of quadratic APN functions via the variable approach}
\label{sec: randomized apn}
Unlike the approach described in Subsection~\ref{sec: bent to apn}, which constructs quadratic $(8,8)$-APN functions by extending suitable quadratic $(8,4)$-bent functions through additional coordinate functions, the present method aims at constructing random quadratic $(8,8)$-APN functions by progressively selecting entries of their look-up tables. Equivalently, this can be viewed as introducing a new input variable to a quadratic $(n,m)$-APN function in order to produce a quadratic $(n+1,m)$-APN function.
Its main advantage is that it does not rely on the working hypothesis from Subsection~\ref{sec: bent to apn}. 
Moreover, since the search is substantially less structured, the functions obtained in this way are expected to provide a more representative sample of the set of quadratic APN functions in dimension~$8$, which is useful for estimating their total number up to EA-equivalence.

The method is based on the following necessary condition: 
If an $(n,n)$-function $F$ is APN, then every function obtained from $F$ by restricting to a lower-dimensional affine subspace of the input space must be APN as well. 
This naturally leads to a recursive construction strategy, which we organize into two phases.
We note that this overall approach is not new, as similar recursive methods have been successfully employed in prior works for constructing 8-bit APN functions (see, e.g., \cite{DBLP:journals/dcc/BeierleLP22}).

\subsection*{Phase I} We first classify quadratic APN $(n,8)$-functions for $n \le 6$, up to EA-equivalence. 
The classification relies on the following result.

\begin{theorem}[Prop.\@ 2.11 of~\cite{DBLP:journals/corr/abs-2501-03922}]
    \label{Thm:dimension_extension}
    Let $F$ be a quadratic APN $(n,m)$-function and $L$ be a linear $(n,m)$-function. 
    Then,
    $$ G(x,x_n) = F(x) + x_n L(x) \,, \quad \text{for } x \in \F_2^n \ \text{and} \  x_n \in \F_2 $$
    is a quadratic $(n+1,m)$-APN function if and only if, for every non-zero $\alpha \in \F_2^n$, we have
    $$ L(\alpha) \neq F(x + \alpha) + F(x) + F(\alpha) + F(0) \,, \quad \forall x \in \F_2^n.$$
\end{theorem}

To construct all quadratic $(n+1,m)$-APN functions using Theorem~\ref{Thm:dimension_extension}, one may proceed as follows. For every quadratic $(n,m)$-APN function $F$, consider all linear $(n,m)$-functions $L$ and test whether the condition of Theorem~\ref{Thm:dimension_extension} is satisfied. Each such valid choice of $L$ yields a quadratic $(n+1,m)$-APN function.

Since our goal is to classify functions only up to EA-equivalence, it suffices to consider a single representative from each EA-equivalence class of quadratic $(n,m)$-APN functions. Indeed, if $F_1$ and $F_2$ are EA-equivalent, then any $(n+1,m)$-APN function $G_1$ obtained from $F_1$ and a suitable linear function $L_1$ will be EA-equivalent to an $(n+1,m)$-APN function $G_2$ obtained from $F_2$ and a unique suitable linear function $L_2$.

\begin{definition}
    \label{def:representative}
    Let $F$ be an $(n,m)$-function.
    The \emph{look-up table} of $F$ is the sequence  
    $$\big(F(0,0,\ldots,0), F(1,0,\ldots,0), F(0,1,\ldots,0),  F(1,1,\ldots,0), \ldots, F(1,1,\ldots,1)\big).$$
    We define the \emph{canonical representative} of an EA-equivalence class $[F]$ as the function $R \in [F]$ whose look-up table $$\big(R(0,0,\ldots,0), R(1,0,\ldots,0), R(0,1,\ldots,0), R(1,1,\ldots,0), \ldots, R(1,1,\ldots,1)\big)$$ is lexicographically smallest among all members of $[F]$.
\end{definition}

To construct all quadratic $(n+1,m)$-APN functions, we iterate over each canonical representative $F$ of the known EA-equivalence classes of quadratic $(n,m)$-APN functions.
For each such canonical representative $F$, we determine all linear $(n,m)$-functions $L$ that satisfy the condition given in Theorem~\ref{Thm:dimension_extension}. 
Each valid choice of $L$ yields a quadratic $(n+1,m)$-APN function $G$.

Note that by taking $F$ as such a canonical representative, we are ensured that we can identify all such canonical representatives $G$ by the following straightforward lemma.

\begin{lemma}\label{lemma:concat}
    Let $G\colon \F_2^n \times \F_2 \to \F_2^m$ be defined by $G(x, x_n) = F(x) + x_n H(x)$, where $H$ is not necessarily linear.
    If $G$ is a canonical representative of $[G]$, then $F$ is a canonical representative of $[F]$.
\end{lemma}
\begin{proof}
    The look-up table of $G$ is the concatenation of look-up tables of $F$ and $F + H$. 
    The fact that $G$ is lexicographically smallest in $[G]$ implies that $F$ is lexicographically smallest in $[F]$.
\end{proof}

\begin{remark}
    Lemma~\ref{lemma:concat} establishes that every canonical representative of a quadratic $(n+1,m)$-APN class must be an extension of a canonical representative of an $(n,m)$-APN class.
    Consequently, by generating all valid extensions from the complete set of canonical $(n,m)$-APN representatives, we are guaranteed that the resulting set of extensions contains the canonical representative for every quadratic $(n+1,m)$-APN equivalence class.
    However, while every generated function $G$ is APN, not every such $G$ is a canonical representative.
    Given that the number of generated APN functions grows rapidly with increasing values of $n$, removing functions that are not canonical representatives is essential for the efficiency of subsequent construction steps, such as building $(n+2, m)$-functions.
    To optimize memory usage, we immediately verify whether each newly constructed function $G$ is a canonical representative, and we store it only if it satisfies this criterion.
\end{remark}

To verify whether a given quadratic $(n+1,m)$-function $G$ is a canonical representative, we systematically search for a lexicographically smaller EA-equivalent function $G' = B \circ (G \circ A + C)$.
Note that for any function $H$, the equivalent function $H' = H + H(0)$ ensures $H'(0) = 0$, and similarly, adding the linear term $x_i H(e_i)$ ensures $H'(e_i) = 0$ for any basis vector $e_i$, which always yields a lexicographically smaller or equal function.
Consequently, if $G$ is a canonical representative, it must satisfy $G(0) = 0$ and $G(e_i) = 0$ for all $0 \le i \le n$.
Besides, to determine if $G$ is a canonical representative, it is sufficient to compare it only against equivalent functions $G'$ that also satisfy $G'(0) = 0$ and $G'(e_i) = 0$ for all $0 \le i \le n$.

Moreover, since $G$ is a quadratic function, any affine input transformation $A(x) + a$ satisfies $G(A(x) + a) = G(A(x)) + D_{a} G(A(x))$, where the derivative $D_{a} G \colon x \mapsto G(x+a) + G(x)$ is an affine function.
This implies that any constant translation in the input can be extracted and absorbed into the outer affine function $C$, allowing us to safely restrict $A$ to be a linear bijection.
Additionally, any constant addition in $C$ can be absorbed into the constant part of the outer affine bijection $B$, meaning we can also restrict $C$ to be strictly linear.
Since both $A$ and $C$ are linear, and we require $G'(0) = 0$, the outer bijection $B$ must map $0$ to $0$, restricting $B$ to a linear bijection as well.
With $A$, $B$, and $C$ restricted to linear mappings, the requirement that $G'(e_i) = 0$ translates to $B(G(A(e_i)) + C(e_i)) = 0$ for all $0 \le i \le n$.
This equality holds if and only if $C(e_i) = G(A(e_i))$ for all $0 \le i \le n$, uniquely determining the linear mapping $C$ for any chosen $A$.

This reduces our search space to iterating over all linear bijections $A$ and all linear bijections $B$.
However, once $A$ is selected and $C$ is subsequently fixed, we construct the output linear bijection $B$ greedily.
Specifically, we sequentially evaluate the internal mapping $G \circ A + C$ and map its novel non-zero outputs to the lexicographically smallest available elements in $\F_2^m$ via $B$.
This greedy construction forces the look-up table of $G'$ to be as lexicographically small as possible.
If at any step this partial table becomes lexicographically strictly smaller than $G$, we immediately abort the entire search and conclude that $G$ is not a canonical representative.
Conversely, if the partial table becomes lexicographically strictly larger than $G$, we prune the search tree, discard the current branch of $A$, and proceed to the next valid partial assignment.

Furthermore, we do not need to determine the linear bijection $A$ completely before checking the lexicographical condition.
Instead, we construct $A$ incrementally by assigning the images of the basis vectors $A(e_0), A(e_1), \dots, A(e_n)$ one at a time.

At each step $k$, assigning $A(e_k)$ fully determines the outputs of the internal mapping $G \circ A + C$ over the subspace spanned by the first $k+1$ basis vectors.
Consequently, this partial information allows us to incrementally define the greedy bijection $B$ and evaluate the partial look-up table of $G'$ up to the elements of this subspace.
If at any intermediate step this partially determined function becomes lexicographically strictly larger than $G$, we instantly prune the search tree.
This early abort mechanism allows us to discard the current partial assignment of $A$ without ever needing to evaluate or guess the remaining basis mappings $A(e_{k+1}), \dots, A(e_n)$.

Since the canonical representative check is the most computationally demanding part of the approach, we aim to avoid constructing redundant EA-equivalent functions.
We therefore apply the following invariance property to reduce the search space for the linear function $L$.

\begin{lemma}
	Let $G$ be a quadratic $(n+1,m)$-function defined by $G(x, x_n) = F(x) + x_n L(x)$, where $F$ is a quadratic $(n,m)$-function and $L$ is a linear $(n,m)$-function.
	For any $\alpha \in \F_2^n$, the extension $G_{\alpha}(x, x_n) = F(x) + x_n L_{\alpha}(x)$ with
	$$ L_{\alpha}(x) = L(x) + F(x+\alpha) + F(x) + F(\alpha) + F(0) $$
	is EA-equivalent to $G(x, x_n)$.
\end{lemma}

\begin{proof}
    Consider the transformation $B\colon \F_2^n \times \F_2 \to \F_2^n \times \F_2$ defined by $B(x, x_n) = (x + x_n \alpha, x_n)$, which is a linear permutation of $\F_2^n \times \F_2$.
    The composition $G \circ B$ is given by
    $$ G(x + x_n \alpha, x_n) = F(x + x_n \alpha) + x_n L(x + x_n \alpha). $$
    Using $F(x + x_n \alpha) = F(x) + x_n\big(F(x+\alpha) + F(x)\big)$ and $L(x + x_n \alpha) = L(x) + x_n L(\alpha)$, we obtain
    \begin{align*}
        G(x + x_n \alpha, x_n) 
        &= F(x) + x_n\big(F(x+\alpha) + F(x) + L(x)\big) + x_n L(\alpha) \\
        &= F(x) + x_n\big(L_\alpha(x) + F(\alpha) + F(0) + L(\alpha)\big) \\
        &= G_\alpha(x, x_n) + x_n\big(F(\alpha) + F(0) + L(\alpha)\big).
    \end{align*}
    Since $G \circ B$ differs from $G_\alpha$ by the linear function $x_n(F(\alpha) + F(0) + L(\alpha))$, it follows that $G$ and $G_\alpha$ are EA-equivalent.
\end{proof}

\begin{remark}
    Note that since $F$ is quadratic, $F(x+\alpha) + F(x) + F(\alpha) + F(0)$ is indeed a linear function in $x$.
    Consequently, for a fixed $F$ and a fixed $L$, all $2^n$ linear functions of the form
    $$ L_\alpha(x) = L(x) + F(x+\alpha) + F(x) + F(\alpha) + F(0) $$
    yield EA-equivalent functions $G_\alpha$.
    Hence, these functions belong to the same EA-equivalence class and need not be considered separately.
    To optimize our search, we only retain the single $L_\alpha$ that produces the lexicographically smallest $G_\alpha$.
\end{remark}

This observation reduces the search space for the linear function $L$ by a factor of $2^n$.
When $F$ is a canonical representative, this reduction in the search space translates into specific constraints on the coefficients of $L$.
Specifically, the number of independent possibilities for the first basis vector of $L$, denoted by $L(1,0, \dots, 0)$, is reduced by a factor of $2^{n-1}$, and the number of possibilities for the second basis vector, $L(0,1, \dots, 0)$, is reduced by a factor of $2$.
This means that $n-1$ coordinates of the first and one coordinate of the second basis vectors of $L$ are fixed to $0$.

Moreover, instead of selecting a fully defined linear function $L$ and subsequently verifying whether it satisfies the conditions of Theorem~\ref{Thm:dimension_extension}, we employ a backtracking algorithm to construct $L$ incrementally.
In this approach, the values of $L$ on the basis vectors of $\F_2^n$ are assigned step by step.
At each step, we continuously verify that the partial assignment does not violate the underlying requirements for the extension $G$ to be APN.
This strategy allows us to prune the search tree significantly, as a single failed partial assignment eliminates an entire branch of candidate linear functions.

To be more explicit, we first enumerate all possible values for $L(e_1)$ that satisfy the necessary condition, with $e_i$ denoting the basis vector of $\F_2^n$ with only $i$-th element being $1$ for $1 \le i \le n$.
For each such choice of $L(e_1)$, we proceed to enumerate all values for $L(e_2)$ that satisfy the conditions for both $L(e_2)$ and $L(e_1 + e_2)$.
This process continues recursively, where at each step $k$, we assign a value to $L(e_k)$ and verify the $2^{k-1}$ conditions for all linear combinations of the basis vectors assigned thus far.
In this manner, the algorithm systematically explores the space of linear $(n,m)$-functions while ensuring that every branch of the search tree remains consistent with the requirements of Theorem~\ref{Thm:dimension_extension}.

It is noteworthy that this recursive construction of the linear function $L$ follows a strategy similar to the one presented in \cite{DBLP:journals/dcc/BeierleLP22}.
We summarize the overall procedure in Algorithm~\ref{alg:APN2}.

\begin{algorithm}
	\caption{Classifying quadratic $(n+1,m)$-APN functions by extending the quadratic $(n,m)$-APN functions}
	\begin{algorithmic}[1]
        \Require Set $\mathcal{Q}_n$ of canonical representative quadratic $(n,m)$-APN functions.
        \Ensure Set $\mathcal{Q}_{n+1}$ of canonical representative quadratic $(n+1,m)$-APN functions.
		\State $\mathcal{Q}_{n+1} \gets \emptyset$.

      	\For{each $F \in \mathcal{Q}_n$}
            \State Initialize $I_0$ with the $2^n \times 2^n$ all-zeros array.
            \State Set $I_0 [\alpha][F(x+\alpha) + F(x) + F(\alpha) + F(0)] \gets 1$ for all $\alpha$ and $x \in \F_2^n$.
            \State Find the $n-1$ conditions on $L(e_1)$ and the single condition on $L(e_2)$.
          	\For{each $b_1 \in \F_2^m$ fulfilling its conditions and with $I_1[e_1][b_1] = 0$}
                \State Set $I_2$ with $I_2[\alpha][b] \gets I_1[\alpha][b] \vee I_1[\alpha + e_1][b + b_1]$ for all $\alpha$ and $b$.
                \For{each $b_2 \in \F_2^m$ fulfilling its condition and with $I_2[e_2][b_2] = 0$}
                    \State Set $I_3$ with $I_3[\alpha][b] \gets I_2[\alpha][b] \vee I_2[\alpha + e_2][b + b_2]$ for all $\alpha$ and $b$.
                    \For{each $b_3 \in \F_2^m$ with $I_3[e_3][b_3] = 0$}
                        \State Set $I_4$ with $I_4[\alpha][b] \gets I_3[\alpha][b] \vee I_3[\alpha + e_3][b + b_3]$ for all $\alpha$ and $b$.
                        \State \ldots
                        \For{each $b_n \in \F_2^m$ with $I_n[e_n][b_n] = 0$}
                            \State Compute $L$ using the basis $(b_1,b_2,\ldots,b_n)$.
                            \State Compute $F'(x,x_n) = F(x) + x_n L(x)$ for all $x$ and $x_n$.
                            \State Check if $F'$ is a canonical representative.
                            \State Add $F'$ to $\mathcal{Q}_{n+1}$ if it is a canonical representative.
                        \EndFor
                        \State \ldots
                    \EndFor
                \EndFor
            \EndFor
		\EndFor
		\State \Return $\mathcal{Q}_{n+1}$.
	\end{algorithmic}
	\label{alg:APN2}
\end{algorithm}

\begin{remark}
    Using the approach given in Algorithm~\ref{alg:APN2}, we obtained the number of EA-equivalence classes of quadratic $(n,m)$-APN functions for the computationally feasible values with $n\leq m$.
    In addition, we reconstructed the known classification results for quadratic $(n,n)$-APN functions in dimensions $n=5,6,7$, obtaining respectively $2$, $13$, and $488$ EA-equivalence classes.
    Successfully reproducing these known results, obtained respectively in \cite{DBLP:journals/dcc/BrinkmannL08,Dillon_Banff:2006,DBLP:journals/iacr/KalginI20a}, serves as strong validation for the correctness of our algorithm.
    Moreover, we classified quadratic $(n,8)$-APN functions for $n\leq 6$ up to EA-equivalence. Table~\ref{table:number_of_classes} provides a comprehensive overview of the number of EA-equivalence classes obtained across computationally feasible dimensions.
    The complete list of functions for this classification is available at \cite{zenodo-classification}.
    \begin{table}[h!]
    	\centering
        \caption{The number of EA-equivalence classes of quadratic $(n,m)$-APN functions}
        \label{table:number_of_classes}
    	\begin{tabular}{c@{\qquad}ccccccc} \toprule
    		\multirow{2}{*}{$m$} & \multicolumn{6}{c}{$n$} \\ \cmidrule{2-7}
    		      & 2 & 3 & 4 &  5 &      6 &     7 \\ \midrule
             2  & 1 &                             \\
             3  & 1 & 1                           \\
             4  & 1 & 1 & 1                       \\
             5  & 1 & 1 & 2 &  2                  \\
             6  & 1 & 1 & 3 &  4 &       13       \\
             7  & 1 & 1 & 3 &  7 &  38\,838 & 488 \\
             8  & 1 & 1 & 3 &  9 & 866\,470 & --- \\
             9  & 1 & 1 & 3 & 10 &      --- & --- \\
            10  & 1 & 1 & 3 & 11 &      --- & --- \\
            11  & 1 & 1 & 3 & 11 &      --- & --- \\ \bottomrule
    	\end{tabular}
    \end{table}
\end{remark}

However, extending this classification to $(7,8)$-functions proved computationally infeasible due to the large search space.
We estimate that there are approximately $2^{48}$ candidate $(7,8)$-APN functions $G$ generated by this approach, making it impossible to perform canonical representative checks for the entire set.
Instead, we employ a randomized construction strategy in the second phase.

\subsection*{Phase II} For each of the $866\,470$ canonical representative quadratic $(6,8)$-APN functions, we perform $2^8$ attempts to construct quadratic $(7,8)$-APN functions, using a procedure similar to that described in the first phase.

In each attempt, we construct the linear function $L$ incrementally.
At each step $k$, a value for $L(e_k)$ is chosen uniformly at random from the set of values that satisfy the APN conditions for all linear combinations of the basis vectors $\{e_1, \dots, e_k\}$ assigned thus far.
If at any step no suitable value exists for $L(e_{k+1})$, the attempt is aborted and a new attempt is started.

If one of the $2^8$ attempts successfully produces a quadratic $(7,8)$-APN function, we then exhaustively search through all suitable linear $(7,8)$-functions in order to construct the corresponding quadratic $(8,8)$-APN function(s).

\begin{remark}
    Using this randomized approach, we generated $92\,955$ quadratic $(8,8)$-APN functions.
    Although this construction method is independent of our working hypothesis, we verified that all obtained functions are indeed extensions of $(8,4)$-bent functions.
    Moreover, this approach alone shows the correctness of two conjectures on the total number of APN functions in dimension eight: that there exist more than $50\,000$ inequivalent quadratic $(8,8)$-APN functions~\cite{YuPerrin2024}, and that their number exceeds the number of inequivalent quadratic $(8,4)$-bent functions, which is $92\,515$, as suggested in~\cite{PolujanPott2022BFA}.
\end{remark}    

\section{Main results}\label{sec:results}
In this section, we summarize the properties of the constructed APN functions, particularly those related to their equivalence to permutations, rank distributions of the non-zero components, and automorphism groups. Additionally, we use the $(8,8)$-APN functions obtained in the previous section to provide an estimate of the total number of APN functions in dimension eight.

\subsection{Constructed APN functions and their properties}\label{sec:constructed_apns}
Using implementations of the approaches introduced in Sections~\ref{sec:millions} and~\ref{sec: randomized apn}, we obtain the main result of this paper:

\begin{theorem}
    There exist at least $3\,808\,491$ CCZ-inequivalent quadratic $(8,8)$-APN functions.
\end{theorem}

In other words, we found $3\,775\,599$ new CCZ-equivalence classes of quadratic APN functions.\footnote{The list can be found here~\cite{dataset}.} See also the data from the project's web page~\cite{apn_project}. We used ortho-derivatives to establish the pairwise CCZ-inequivalence of our found functions.  Notably, the majority of these functions (millions) were obtained using the approach described in Section~\ref{sec: bent to apn}. None of our found functions is CCZ-equivalent to a permutation. For establishing the inequivalence using ortho-derivatives and verifying the inequivalence to a permutation, we used the \texttt{sboxU} tool~\cite{Perrin_sboxu}. The latter verification is based on the fact that the CCZ-equivalence of an $(n,n)$-function $F$ to a permutation is characterized by the existence of trivially intersecting $n$-dimensional vector spaces in the Walsh zeroes of $F$~\cite{perrin_equivalence_permutation,DBLP:journals/ccds/ChaseL22}.  

The \emph{rank distribution} of a quadratic $(n,n)$-APN function $F$ is the multiset counting the ranks of all non-zero component functions $f \in \comp(F)\setminus{0}$. For even $n$, rank distributions are completely determined only for $n \le 6$, where the full classification of quadratic APN functions is known. For $n \ge 8$, determining the possible rank distributions of quadratic $(n,n)$-APN functions remains a difficult problem. In a recent work~\cite{Koelsch2025}, the possible rank distributions of quadratic $(n,n)$-APN functions were determined in the following theorem.

\begin{theorem}\cite[{Theorem IV.5.}]{Koelsch2025}\label{th:rank_distributions}
    Let $F$ be a quadratic $(8,8)$-APN function. Then, $F$ has one of the following rank distributions:
\begin{itemize}
    \item $\{* \ 8^{190}, 6^{64}, 4^0, 2^{1} \ *\}$,
    \item $\{* \ 8^{238-4i}, 6^{5i}, 4^{17-i}, 2^0 \ *\}$, for $1\leq i\leq 17$.
\end{itemize}
\end{theorem}

In Table~\ref{table:ranks_combined}, we summarize the rank distributions (by listing the multiplicities of the quadratic components with ranks $r=8,6,4,2$) of the already known $32\,892$ quadratic $(8,8)$-APN functions and compare them with those of the newly constructed $3\,775\,599$ examples.

\begin{table}[h!]
\centering
\caption{Rank distributions of quadratic $(8,8)$-APN functions}
\begin{subtable}{0.48\textwidth}
\centering
\caption{Known $32\,892$ APN functions}
\begin{tabular}{c c c c c r}
\hline
Type & $8$ & $6$ & $4$ & $2$ & Count \\
\hline
$1$ & $170$ & $85$ & $0$ & $0$ & $27\,442$ \\
$2$ & $174$ & $80$ & $1$ & $0$ & $4\,913$ \\
$3$ & $178$ & $75$ & $2$ & $0$ & $361$ \\
$4$ & $182$ & $70$ & $3$ & $0$ & $146$ \\
$5$ & $186$ & $65$ & $4$ & $0$ & $26$ \\
$6$ & $190$ & $64$ & $0$ & $1$ & $4$ \\
\hline
\end{tabular}
\label{table:ranks_known_before}
\end{subtable}
\hfill
\begin{subtable}{0.48\textwidth}
\centering
\caption{New $3\,775\,599$ APN functions}
\begin{tabular}{c c c c c r}
\hline
Type & $8$ & $6$ & $4$ & $2$ & Count \\
\hline
$1$ & $170$ & $85$ & $0$ & $0$ & $3\,454\,490$ \\
$2$ & $174$ & $80$ & $1$ & $0$ & $313\,910$ \\
$3$ & $178$ & $75$ & $2$ & $0$ & $7\,136$ \\
$4$ & $182$ & $70$ & $3$ & $0$ & $63$ \\
$5$ & $186$ & $65$ & $4$ & $0$ & $0$ \\
$6$ & $190$ & $64$ & $0$ & $1$ & $0$ \\
\hline
\end{tabular}
\label{table:ranks_new}
\end{subtable}
\label{table:ranks_combined}
\end{table}

As one can see from Tables~\ref{table:ranks_known_before} and~\ref{table:ranks_new}, the vast majority of the newly constructed functions has the classical rank distribution corresponding to the case with the $2/3$ of nonzero bent components and $1/3$ of semi-bent components. No new ``rare'' APN functions were found with our approach.\footnote{There are exactly $4$ CCZ-equivalence classes of quadratic $(8,8)$-APN functions with $190$ bent components~\cite{DBLP:journals/dcc/BeierleLP22}.} Moreover, no rank distributions of the form $\{* \ 8^{238-4i}, 6^{5i}, 4^{17-i}, 2^0 \ *\}$ for $1\leq i\leq 12$ have been found. Thus, Open Problem VI.1 in~\cite{Koelsch2025} concerning the existence of such functions remains unresolved.

Another important parameter of a quadratic $(n,n)$-APN function $F$  is its \emph{automorphism group} $\operatorname{Aut}(F)$, defined as the group of all $\F_2$-affine permutations $\sigma$ from $\F_2^n\times \F_2^n$ to itself such that the graph  $\mathcal{G}_F=\{ (x,F(x))\colon x \in \F_2^n\}$ is invariant under $\sigma$. It is well-known that for CCZ-equivalent functions, the automorphism groups are conjugate, so their orders are the same. The automorphism group of a quadratic $(n,n)$-function $F$ (including APN) contains an elementary
abelian group of order $2^n$ (see~\cite[Theorem 1]{EdelP09} for the exact description of its elements). If $|\operatorname{Aut}(F)|=2^n$, we say that the automorphism group of $F$ is \emph{trivial}.
 
In Table~\ref{table:aut_combined}, we summarize orders of the automorphism groups of quadratic $(8,8)$-APN functions of the already known $32\,892$ quadratic $(8,8)$-APN functions and compare them with those of the newly constructed $3\,775\,599$ examples. 

\begin{table}[h!]
\centering
\caption{Orders of the automorphism groups of quadratic $(8,8)$-APN functions}
\begin{subtable}{\textwidth}
\centering
\caption{Known $32\,892$ APN functions}
\begin{tabular}{cccccccccc}
\hline
$|\operatorname{Aut}(F)|$
& $2^8$
& $2^9$
& $2^8 \cdot 3$
& $2^{10}$
& $2^8 \cdot 5$
& $2^9 \cdot 3$
& $2^8 \cdot 7$
& $2^{11}$ 
& $2^8 \cdot 3^2$ \\
\hline
Count
& $20\,118$
& $69$
& $12\,647$
& $8$
& $3$
& $11$
& $9$
& $3$ 
& $3$ \\
\hline
\end{tabular}
\begin{tabular}{ccccccc}
\hline
$|\operatorname{Aut}(F)|$
& $2^{10} \cdot 3$
& $2^8 \cdot 3 \cdot 5$
& $2^{11} \cdot 3$
& $2^{10} \cdot 3 \cdot 5$
& $2^{10} \cdot 3^2 \cdot 5$
& $2^{11} \cdot 3 \cdot 5 \cdot 17$ \\
\hline
Count
& $9$
& $5$
& $2$
& $1$
& $2$
& $2$ \\
\hline
\end{tabular}
\vspace{2mm}
\label{table:aut_known_before}
\end{subtable}
\begin{subtable}{\textwidth}
\centering
\caption{New $3\,775\,599$ APN functions}
\begin{tabular}{ccc}
\hline
$|\operatorname{Aut}(F)|$ & $2^8$ & $2^8 \cdot 3$ \\
\hline
Count & $3\,775\,598$ & $1$ \\
\hline
\end{tabular}
\label{table:aut_new}
\end{subtable}
\label{table:aut_combined}
\end{table}

To compute the automorphism groups of the known and newly constructed functions, we implemented in \texttt{Magma}~\cite{MAGMA:1997} the coding-theoretic approach described in~\cite{EdelP09Equivalence}. Remarkably, all but one function in our list have a trivial automorphism group. The only function with a non-trivial automorphism group is function no.~$748\,877$ in our list, whose automorphism group has size $2^8 \cdot 3$.

\subsection{On the total number of inequivalent quadratic APN functions in dimension eight}\label{sec:total_number}
Finally, we would like to estimate the total number of inequivalent quadratic $(8,8)$-APN functions. To do so, we assume that the $92\,955$ functions generated by the method described in Section~\ref{sec: randomized apn} constitutes a uniformly random sample (chosen with replacement) of all CCZ-equivalence classes of quadratic APN functions in 8 variables. Let us denote this latter quantity by $N$. Estimating $N$ then boils down to an \emph{inverse coupon collector's problem}: Suppose we have a uniform sample of $t$ objects (here APN functions) from $N$ objects, chosen with replacement. Suppose we have $\ell$ distinct objects (here, distinct ortho-derivative labels to indicate CCZ-inequivalent functions) in our sample. As explained in~\cite{30d68731-14d3-304e-86d4-94c9bf3957ff}, the maximum-likelihood estimator for $N$ is given by $$\argmax_{N > {\ell-1}}\{\binom{N}{\ell}/N^t\},$$ and the function $\mathbb{R} \rightarrow \mathbb{R}, N \mapsto \binom{N}{\ell}/N^t$ has a unique local maximum for $N > \ell-1$. In our case, we have $(t,\ell) = (92\,955,92\,253)$ and a local maximum for $N$ around $N = 6\,123\,206$.

To obtain another estimate for $N$, we checked how many of the $t$ functions belong to one of the $M = 3\,776\,451
$ equivalence classes known before or found with the approach explained in Section~\ref{sec: bent to apn}. In our case, $t' = 32\,286$ out of the $t$ functions do \emph{not} belong to any of those $M$ equivalence classes. The expected value for $t'$ is given by $(1-M/N)t$, so an estimate for $N$ is calculated as $N = tM/(t-t') \approx 5\,786\,151$.

\section{Conclusion and open problems}\label{sec:concl_op}
We have constructed the by far largest known number of quadratic $(8,8)$-APN functions, which, according to our estimation, is expected to be more than half of the total number. To conclude the paper, we would like to give the following list of problems, which would hopefully further stimulate future investigation of APN functions.
\begin{enumerate}
	\item Further investigate whether every quadratic $(n,n)$-APN function, with $n$ even, is an extension of a quadratic vectorial bent $(n,n/2)$-function. Obtaining tighter bounds than those by Clark and Suen~\cite{clark1999probability} on the probability that a subset of $\F_2^n$ of given size contains a vector space of dimension $k$ would be a very interesting problem.
	\item While millions of new examples were constructed, our estimate indicates that we are still far from what seems to be the total number of quadratic $(n,n)$-APN functions. We believe it is still interesting to develop further efficient construction methods of quadratic $(n,n)$-APN functions in dimensions $n \ge 8$, aiming to classify the first open case $n=8$. This problem seems difficult to solve at the moment due to a lack of understanding of the combinatorial structure of quadratic APN functions (e.g., possible rank distributions of component functions and relations between them, as well as possible automorphism groups). Additionally, it is important to develop further fast-to-compute and precise invariants for $(n,m)$-functions, since, as indicated by this work, early application of such invariants helps to reduce the combinatorial explosion early.
	\item This work provides further evidence that $(k,n)$-APN functions with $k<n$ can serve as building blocks for constructing $(n,n)$-APN functions. While $(n,n)$-APN functions are well studied, $(k,n)$-APN functions remain less understood. Developing construction methods and studying their properties is therefore an important direction for future research.
\end{enumerate}

\section*{Acknowledgments}
The results given in this article were presented at the \emph{10th International Workshop on Boolean Functions and their Applications (BFA), September 1-5, 2025, Larnaca, Cyprus}.

Philippe Langevin is partially supported by the \emph{French Agence Nationale de la Recherche} through the SWAP project under the Contract ANR-21-CE39-0012. Shahram Rasoolzadeh is funded by the \emph{ERC} project 101097056 (SYMTRUST). Christof Beierle and Alexandr Polujan are supported by the Deutsche Forschungsgemeinschaft (DFG, German Research Foundation) --- project number 581467925.

The authors would like to thank Claude Carlet, Faruk G{\"{o}}loglu and Lukas K\"olsch for useful information.


\begin{thebibliography}{10}
\providecommand{\url}[1]{#1}
\csname url@samestyle\endcsname
\providecommand{\newblock}{\relax}
\providecommand{\bibinfo}[2]{#2}
\providecommand{\BIBentrySTDinterwordspacing}{\spaceskip=0pt\relax}
\providecommand{\BIBentryALTinterwordstretchfactor}{4}
\providecommand{\BIBentryALTinterwordspacing}{\spaceskip=\fontdimen2\font plus
\BIBentryALTinterwordstretchfactor\fontdimen3\font minus \fontdimen4\font\relax}
\providecommand{\BIBforeignlanguage}[2]{{%
\expandafter\ifx\csname l@#1\endcsname\relax
\typeout{** WARNING: IEEEtranS.bst: No hyphenation pattern has been}%
\typeout{** loaded for the language `#1'. Using the pattern for}%
\typeout{** the default language instead.}%
\else
\language=\csname l@#1\endcsname
\fi
#2}}
\providecommand{\BIBdecl}{\relax}
\BIBdecl

\bibitem{Arshad_PhD_Thesis}
R.~Arshad, ``Contributions to the theory of almost perfect nonlinear functions,'' Ph.D. dissertation, Otto-von-Guericke-Universität Magdeburg, Fakultät für Mathematik, 2018.

\bibitem{DBLP:journals/ffa/BeierleCLP22}
C.~Beierle, C.~Carlet, G.~Leander, and L.~Perrin, ``A further study of quadratic {APN} permutations in dimension nine,'' \emph{Finite Fields Their Appl.}, vol.~81, p. 102049, 2022.

\bibitem{apn_project}
\BIBentryALTinterwordspacing
C.~Beierle, P.~Langevin, G.~Leander, A.~Polujan, and S.~Rasoolzadeh, ``Quadratic {APN} functions in dimension 8.'' [Online]. Available: \url{https://langevin.univ-tln.fr/data/apns/}
\BIBentrySTDinterwordspacing

\bibitem{dataset}
\BIBentryALTinterwordspacing
C.~Beierle, P.~Langevin, G.~Leander, A.~Polujan, and S.~Rasoolzadeh, ``Millions of inequivalent quadratic {APN} functions in eight variables,'' Zenodo, 2025. [Online]. Available: \url{https://doi.org/10.5281/zenodo.16752428}
\BIBentrySTDinterwordspacing

\bibitem{zenodo-classification}
\BIBentryALTinterwordspacing
C.~Beierle, P.~Langevin, G.~Leander, A.~Polujan, and S.~Rasoolzadeh, ``Classification of quadratic {APN} functions up to {EA}-equivalence,'' 2026. [Online]. Available: \url{https://doi.org/10.5281/zenodo.22010376}
\BIBentrySTDinterwordspacing

\bibitem{DBLP:journals/corr/abs-2009-07204}
C.~Beierle and G.~Leander, ``New instances of quadratic {APN} functions,'' \emph{{IEEE} Trans. Inf. Theory}, vol.~68, no.~1, pp. 670--678, 2022.

\bibitem{DBLP:journals/dcc/BeierleLP22}
C.~Beierle, G.~Leander, and L.~Perrin, ``Trims and extensions of quadratic {APN} functions,'' \emph{Des. Codes Cryptogr.}, vol.~90, no.~4, pp. 1009--1036, 2022.

\bibitem{berger2006almost}
T.~P. Berger, A.~Canteaut, P.~Charpin, and Y.~Laigle-Chapuy, ``On almost perfect nonlinear functions over $\mathbb{F}_2^n$,'' \emph{IEEE Transactions on Information Theory}, vol.~52, no.~9, pp. 4160--4170, 2006.

\bibitem{BethDing1994}
T.~Beth and C.~Ding, ``On almost perfect nonlinear permutations,'' in \emph{Advances in cryptology -- EUROCRYPT `93. Proceedings of the workshop on the theory and application of cryptographic techniques, Lofthus, Norway, May 23--27, 1993. Proceedings}.\hskip 1em plus 0.5em minus 0.4em\relax Berlin: Springer, 1994, pp. 65--76.

\bibitem{DBLP:conf/asiacrypt/BonnetainPT19}
X.~Bonnetain, L.~Perrin, and S.~Tian, ``Anomalies and vector space search: Tools for s-box analysis,'' in \emph{Advances in Cryptology - {ASIACRYPT} 2019 - 25th International Conference on the Theory and Application of Cryptology and Information Security, Kobe, Japan, December 8-12, 2019, Proceedings, Part {I}}, ser. Lecture Notes in Computer Science, S.~D. Galbraith and S.~Moriai, Eds., vol. 11921.\hskip 1em plus 0.5em minus 0.4em\relax Springer, 2019, pp. 196--223.

\bibitem{MAGMA:1997}
W.~Bosma, J.~Cannon, and C.~Playoust, ``The {M}agma algebra system. {I}. {T}he user language,'' \emph{Journal of Symbolic Computation}, vol.~24, no. 3-4, pp. 235--265, 1997.

\bibitem{DBLP:journals/dcc/BrinkmannL08}
M.~Brinkmann and G.~Leander, ``On the classification of {APN} functions up to dimension five,'' \emph{Des. Codes Cryptogr.}, vol.~49, no. 1-3, pp. 273--288, 2008.

\bibitem{Dillon-perm}
K.~A. Browning, J.~F. Dillon, M.~T. McQuistan, and A.~J. Wolfe, ``\BIBforeignlanguage{English}{An {APN} permutation in dimension six},'' in \emph{\BIBforeignlanguage{English}{Finite fields. Theory and applications. Proceedings of the 9th international conference on finite fields and applications, Dublin, Ireland, July 13--17, 2009}}.\hskip 1em plus 0.5em minus 0.4em\relax Providence, RI: American Mathematical Society (AMS), 2010, pp. 33--42.

\bibitem{BudaghyanC2010}
L.~Budaghyan and C.~Carlet, ``\BIBforeignlanguage{English}{{CCZ}-equivalence of single and multi output {Boolean} functions},'' in \emph{\BIBforeignlanguage{English}{Finite fields. Theory and applications. Proceedings of the 9th international conference on finite fields and applications, Dublin, Ireland, July 13--17, 2009}}.\hskip 1em plus 0.5em minus 0.4em\relax Providence, RI: American Mathematical Society (AMS), 2010, pp. 43--54.

\bibitem{Koelsch2025}
S.~H. Bénéteau, N.~Goluboff, L.~Kölsch, and D.~Vaghasiya, ``On the {W}alsh spectra of quadratic {APN} functions,'' \emph{IEEE Transactions on Information Theory}, vol.~72, no.~7, pp. 5207--5216, 2026.

\bibitem{CanteautCP2022}
A.~Canteaut, A.~Couvreur, and L.~Perrin, ``Recovering or testing extended-affine equivalence,'' \emph{IEEE Transactions on Information Theory}, vol.~68, no.~9, pp. 6187--6206, 2022.

\bibitem{perrin_equivalence_permutation}
A.~Canteaut and L.~Perrin, ``On ccz-equivalence, extended-affine equivalence and function twisting,'' in \emph{The 14th International Conference on Finite Fields and their Applications}, Vancouver, Canada, 2019.

\bibitem{Carlet11}
C.~Carlet, ``Relating three nonlinearity parameters of vectorial functions and building {APN} functions from bent functions,'' \emph{Des. Codes Cryptogr.}, vol.~59, no. 1-3, pp. 89--109, 2011.

\bibitem{Carlet2021_Book}
C.~Carlet, \emph{{B}oolean {F}unctions for {C}ryptography and {C}oding {T}heory}.\hskip 1em plus 0.5em minus 0.4em\relax Cambridge University Press, 2021.

\bibitem{DBLP:journals/dcc/CarletCZ98}
C.~Carlet, P.~Charpin, and V.~A. Zinoviev, ``Codes, bent functions and permutations suitable for {DES}-like cryptosystems,'' \emph{Des. Codes Cryptogr.}, vol.~15, no.~2, pp. 125--156, 1998.

\bibitem{CarletPicek2021}
C.~Carlet and S.~Picek, ``On the exponents of {APN} power functions and {S}idon sets, sum-free sets, and {D}ickson polynomials,'' \emph{Advances in Mathematics of Communications}, vol.~17, no.~6, pp. 1507--1525, 2023.

\bibitem{DBLP:journals/ccds/ChaseL22}
B.~Chase and P.~Lisonek, ``Construction of {APN} permutations via {Walsh} zero spaces,'' \emph{Cryptogr. Commun.}, vol.~14, no.~6, pp. 1345--1357, 2022.

\bibitem{clark1999probability}
W.~E. Clark and S.~Suen, ``On the probability that a {$t$}-subset of a finite vector space contains an {$r$}-subspace -- with applications to short, light codewords in a {BCH} code,'' \emph{Congr. Numerantium}, vol. 137, pp. 139--159, 1999.

\bibitem{CzerwinskiP26}
I.~Czerwinski and A.~Pott, ``On large {S}idon sets,'' \emph{J. Comb. Theory, Ser. {A}}, vol. 220, p. 106129, 2026.

\bibitem{30d68731-14d3-304e-86d4-94c9bf3957ff}
\BIBentryALTinterwordspacing
B.~Dawkins, ``Siobhan's problem: The coupon collector revisited,'' \emph{The American Statistician}, vol.~45, no.~1, pp. 76--82, 1991. [Online]. Available: \url{http://www.jstor.org/stable/2685247}
\BIBentrySTDinterwordspacing

\bibitem{Dillon1999}
J.~F. Dillon, ``Multiplicative difference sets via additive characters,'' \emph{Designs, Codes and Cryptography}, vol.~17, no.~1, pp. 225--235, 1999.

\bibitem{Dillon_Banff:2006}
J.~F. Dillon, ``Slides from talk given at \textquotedblleft {P}olynomials over finite fields and applications\textquotedblright, held at {B}anff {I}nternational {R}esearch {S}tation,'' 2006.

\bibitem{DillonDobbertin2004}
J.~F. Dillon and H.~Dobbertin, ``New cyclic difference sets with {Singer} parameters,'' \emph{Finite Fields and Their Applications}, vol.~10, no.~3, pp. 342--389, 2004.

\bibitem{DingMT19}
C.~Ding, A.~Munemasa, and V.~D. Tonchev, ``Bent vectorial functions, codes and designs,'' \emph{{IEEE} Trans. Inf. Theory}, vol.~65, no.~11, pp. 7533--7541, 2019.

\bibitem{Dobbertin1999:Niho}
H.~Dobbertin, ``Almost perfect nonlinear power functions on {{\(\mathrm{GF}(2^n)\)}}: the {Niho} case.'' \emph{Inf. Comput.}, vol. 151, no. 1-2, pp. 57--72, 1999.

\bibitem{Dobbertin1999:Welch}
H.~Dobbertin, ``\BIBforeignlanguage{English}{Almost perfect nonlinear power functions on {{\(\mathrm{GF}(2^n)\)}}: the {Welch} case},'' \emph{\BIBforeignlanguage{English}{IEEE Trans. Inf. Theory}}, vol.~45, no.~4, pp. 1271--1275, 1999.

\bibitem{Dobbertin2001:DivBy5}
H.~Dobbertin, ``\BIBforeignlanguage{English}{Almost perfect nonlinear power functions on {{\(\text{GF} (2^n)\)}}: {A} new case for {{\(n\)}} divisible by 5},'' in \emph{\BIBforeignlanguage{English}{Finite fields and applications. Proceedings of the fifth international conference on finite fields and applications \(F_q5\), University of Augsburg, Germany, August 2--6, 1999}}.\hskip 1em plus 0.5em minus 0.4em\relax Berlin: Springer, 2001, pp. 113--121.

\bibitem{DobbertinMMPW2003}
H.~Dobbertin, D.~Mills, E.~N. M{\"{u}}ller, A.~Pott, and W.~Willems, ``{APN} functions in odd characteristic,'' \emph{Discret. Math.}, vol. 267, no. 1-3, pp. 95--112, 2003.

\bibitem{Edel2011}
\BIBentryALTinterwordspacing
Y.~Edel, ``Quadratic {APN} functions as subspaces of alternating bilinear forms,'' in \emph{Proceedings of the Contact Forum Coding Theory and Cryptography III at The Royal Flemish Academy of Belgium for Science and the Arts 2009}, 2011, pp. 11--24. [Online]. Available: \url{http://www.yvesedel.de/Papers/ContactForum09.html}
\BIBentrySTDinterwordspacing

\bibitem{EdelKP2006}
Y.~Edel, G.~M. Kyureghyan, and A.~Pott, ``A new {APN} function which is not equivalent to a power mapping,'' \emph{{IEEE} Trans. Inf. Theory}, vol.~52, no.~2, pp. 744--747, 2006.

\bibitem{EdelP09}
Y.~Edel and A.~Pott, ``A new almost perfect nonlinear function which is not quadratic,'' \emph{Adv. Math. Commun.}, vol.~3, no.~1, pp. 59--81, 2009.

\bibitem{EdelP09Equivalence}
\BIBentryALTinterwordspacing
Y.~Edel and A.~Pott, ``On the equivalence of nonlinear functions,'' in \emph{Enhancing Cryptographic Primitives with Techniques from Error Correcting Codes}, 2009, pp. 87--103. [Online]. Available: \url{http://hdl.handle.net/1854/LU-710528}
\BIBentrySTDinterwordspacing

\bibitem{Gold1968}
R.~Gold, ``Maximal recursive sequences with 3-valued recursive cross-correlation functions,'' \emph{{IEEE} Trans. Inf. Theory}, vol.~14, no.~1, pp. 154--156, 1968.

\bibitem{GologluK2025_CT}
F.~G{\"o}lo{\u{g}}lu and L.~K{\"o}lsch, ``Equivalences of biprojective almost perfect nonlinear functions,'' \emph{Comb. Theory}, vol.~5, no.~3, p.~28, 2025, id/No 7.

\bibitem{GologluL20}
F.~G{\"{o}}loglu and P.~Langevin, ``Almost perfect nonlinear families which are not equivalent to permutations,'' \emph{Finite Fields Their Appl.}, vol.~67, p. 101707, 2020.

\bibitem{GologluPavlu2021}
F.~G{\"o}lo{\u g}lu and J.~Pavlu, ``On {CCZ}-inequivalence of some families of almost perfect nonlinear functions to permutations,'' \emph{Cryptography and Communications}, vol.~13, no.~3, pp. 377--391, 2021.

\bibitem{Hou06}
X.~Hou, ``Affinity of permutations of {$\mathbb{F}_2^n$},'' \emph{Discret. Appl. Math.}, vol. 154, no.~2, pp. 313--325, 2006.

\bibitem{IchanskaBKY2025}
N.~Ichanska, S.~Berg, N.~Kaleyski, and Y.~Yu, ``Pushing the {QAM} method for finding {APN} functions further,'' \emph{Cryptography and Communications}, 2025.

\bibitem{JanwaWilson1993}
H.~Janwa and R.~Wilson, ``Hyperplane sections of fermat varieties in {$\mathbb{P}^3$} in char. 2 and some applications to cyclic codes,'' in \emph{Proceedings of AAECC-10 Conference}, ser. Lecture Notes in Computer Science, vol. 673, 1993, pp. 180--194.

\bibitem{DBLP:journals/ccds/Kaleyski21}
N.~S. Kaleyski, ``Invariants for {EA-} and {CCZ-}equivalence of {APN} and {AB} functions,'' \emph{Cryptogr. Commun.}, vol.~13, no.~6, pp. 995--1023, 2021.

\bibitem{DBLP:journals/iacr/KalginI20a}
K.~Kalgin and V.~Idrisova, ``The classification of quadratic {APN} functions in 7 variables and combinatorial approaches to search for {APN} functions,'' \emph{Cryptogr. Commun.}, vol.~15, no.~2, pp. 239--256, 2023.

\bibitem{Kasami1971}
T.~Kasami, ``The weight enumerators for several classes of subcodes of the 2nd order binary {Reed}-{Muller} codes,'' \emph{Inf. Control}, vol.~18, pp. 369--394, 1971.

\bibitem{KaspersZhou2021}
C.~Kaspers and Y.~Zhou, ``The number of almost perfect nonlinear functions grows exponentially,'' \emph{J. Cryptol.}, vol.~34, no.~1, p.~4, 2021.

\bibitem{KyureghyanP08}
G.~M. Kyureghyan and A.~Pott, ``Some theorems on planar mappings,'' in \emph{Arithmetic of Finite Fields, 2nd International Workshop, {WAIFI} 2008, Siena, Italy, July 6-9, 2008, Proceedings}, ser. Lecture Notes in Computer Science, J.~von~zur Gathen, J.~L. Ima{\~{n}}a, and {\c{C}}.~K. Ko{\c{c}}, Eds.\hskip 1em plus 0.5em minus 0.4em\relax Springer, 2008, pp. 117--122.

\bibitem{langevin}
P.~Langevin, ``Classification of {APN} cubics in dimension 6 over {GF(2)},'' \url{http://langevin.univ-tln.fr/project/apn-6/apn-6.html}, 2012, accessed February 10, 2026.

\bibitem{Li_Kaleyski_2024}
K.~Li and N.~Kaleyski, ``Two new infinite families of {APN} functions in trivariate form,'' \emph{{IEEE} Trans. Inf. Theory}, vol.~70, no.~2, pp. 1436--1452, 2024.

\bibitem{MacWilliamsSloane}
F.~J. MacWilliams and N.~J.~A. Sloane, \emph{The Theory of Error-Correcting Codes. {I}}.\hskip 1em plus 0.5em minus 0.4em\relax Amsterdam: North-Holland Mathematical Library, Vol. 16, North-Holland Publishing Co., 1977.

\bibitem{Nyberg91}
K.~Nyberg, ``Perfect nonlinear {S}-{B}oxes,'' in \emph{Advances in Cryptology - {EUROCRYPT} '91, Workshop on the Theory and Application of of Cryptographic Techniques, Brighton, UK, April 8-11, 1991, Proceedings}, ser. Lecture Notes in Computer Science, D.~W. Davies, Ed., vol. 547.\hskip 1em plus 0.5em minus 0.4em\relax Springer, 1991, pp. 378--386.

\bibitem{Nyberg93}
K.~Nyberg, ``Differentially uniform mappings for cryptography,'' in \emph{Advances in Cryptology - {EUROCRYPT} '93, Workshop on the Theory and Application of Cryptographic Techniques, Lofthus, Norway, May 23-27, 1993, Proceedings}, ser. Lecture Notes in Computer Science, T.~Helleseth, Ed., vol. 765.\hskip 1em plus 0.5em minus 0.4em\relax Springer, 1993, pp. 55--64.

\bibitem{DBLP:conf/fse/Nyberg94}
K.~Nyberg, ``S-boxes and round functions with controllable linearity and differential uniformity,'' in \emph{Fast Software Encryption: Second International Workshop. Leuven, Belgium, 14-16 December 1994, Proceedings}, ser. Lecture Notes in Computer Science, B.~Preneel, Ed., vol. 1008.\hskip 1em plus 0.5em minus 0.4em\relax Springer, 1994, pp. 111--130.

\bibitem{DBLP:conf/crypto/NybergK92}
K.~Nyberg and L.~R. Knudsen, ``Provable security against differential cryptanalysis,'' in \emph{Advances in Cryptology - {CRYPTO} '92}, ser. Lecture Notes in Computer Science, E.~F. Brickell, Ed., vol. 740.\hskip 1em plus 0.5em minus 0.4em\relax Springer, 1992, pp. 566--574.

\bibitem{Perrin_sboxu}
L.~Perrin, ``{sboxU},'' \emph{GitHub repository}, 2025, v1.3.1 availabe via \url{https://github.com/lpp-crypto/sboxU}.

\bibitem{PolujanPott2022BFA}
A.~Polujan and A.~Pott, ``Towards the classification of quadratic vectorial bent functions in 8 variables,'' \url{https://boolean.w.uib.no/bfa-2022/}, 2022, paper 12.

\bibitem{PolujanP21}
A.~A. Polujan and A.~Pott, ``On design-theoretic aspects of boolean and vectorial bent functions,'' \emph{{IEEE} Trans. Inf. Theory}, vol.~67, no.~2, pp. 1027--1037, 2021.

\bibitem{DBLP:journals/dcc/Pott16}
A.~Pott, ``Almost perfect and planar functions,'' \emph{Des. Codes Cryptogr.}, vol.~78, no.~1, pp. 141--195, 2016.

\bibitem{ROTHAUS1976}
O.~S. Rothaus, ``On ``bent'' functions,'' \emph{Journal of Combinatorial Theory, Series A}, vol.~20, no.~3, pp. 300--305, 1976.

\bibitem{DBLP:journals/dcc/Taniguchi19a}
H.~Taniguchi, ``On some quadratic {APN} functions,'' \emph{Des. Codes Cryptogr.}, vol.~87, no.~9, pp. 1973--1983, 2019.

\bibitem{DBLP:journals/corr/abs-2501-03922}
H.~Taniguchi, A.~Polujan, A.~Pott, and R.~Arshad, ``Changing almost perfect nonlinear functions on affine subspaces of small codimensions,'' \emph{The Electronic Journal of Combinatorics}, vol.~32, no.~4, p. \#P4.61, 2025.

\bibitem{weng2013quadratic}
G.~Weng, Y.~Tan, and G.~Gong, ``On quadratic almost perfect nonlinear functions and their related algebraic object,'' in \emph{Workshop on Coding and Cryptography, WCC.}, 2013.

\bibitem{yoshiara2008dimensional}
S.~Yoshiara, ``Dimensional dual hyperovals associated with quadratic {APN} functions,'' \emph{Innovations in Incidence Geometry: Algebraic, Topological and Combinatorial}, vol.~8, no.~1, pp. 147--169, 2008.

\bibitem{Yoshiara10}
S.~Yoshiara, ``Notes on {APN} functions, semibiplanes and dimensional dual hyperovals,'' \emph{Des. Codes Cryptogr.}, vol.~56, no. 2-3, pp. 197--218, 2010.

\bibitem{YuPerrin2024}
Y.~Yu and L.~Perrin, ``Constructing more quadratic {APN} functions with the {QAM} method,'' \emph{Cryptography and Communications}, vol.~14, no.~6, pp. 1359--1369, 2022.

\bibitem{DBLP:journals/dcc/YuWL14}
Y.~Yu, M.~Wang, and Y.~Li, ``A matrix approach for constructing quadratic {APN} functions,'' \emph{Des. Codes Cryptogr.}, vol.~73, no.~2, pp. 587--600, 2014.

\end{thebibliography}


\end{document}